\documentclass[9pt,a4paper]{article}

\newcommand{\R}{\mathbb{R}}
\newcommand{\M}{\mathcal{M}}
\newcommand{\T}{\mathcal{T}}

\newcommand{\f}{\mathrm{f}}

\usepackage{tikz}
\usepackage{pgfplots} 
\usepackage{pgfgantt}
\pgfplotsset{compat=newest}
\pgfplotsset{compat=newest} 
\pgfplotsset{plot coordinates/math parser=false}
\usepackage{pgf}
\usepackage{cite}
\usepackage{amsmath}
\usepackage{cleveref}
\usepackage[markup=nocolor]{changes}
\usepackage{xcolor}
\usepackage{lipsum}
\usepackage{amsfonts}
\usepackage{graphicx}
\usepackage{epstopdf}
\usepackage{algorithmic}

\usepackage{amsthm}
\usepackage{fancyhdr}

\newtheorem{definition}{Definition}
\newtheorem{remark}{Remark}
\newtheorem{corollary}{Corollary}
\newtheorem{proposition}{Proposition}
\newtheorem{theorem}{Theorem}
\newtheorem{lemma}{Lemma}

\usepackage{pgfplots}
\usepgfplotslibrary{groupplots}

\title{On Differential Geometric Formulations of Slow Invariant Manifold Computation: Geodesic Stretching and Flow Curvature\thanks{Submitted to the editors DATE.
		This work supported by funding of the Klaus Tschira Foundation (project 00.003.2019).}}

\author{Dirk Lebiedz\thanks{Institute of Numerical Mathematics, Ulm, Germany
		(dirk.lebiedz@uni-ulm.de), (johannespoppe92@gmail.com)}
			\and Johannes Poppe\footnotemark[2]}

\providecommand{\keywords}[1]{\textbf{Key words ---} #1}
\providecommand{\AMS}[1]{\textbf{AMS subject classifications ---} #1}

\date{\today}

\renewcommand{\d}{\mathrm{d}}

\begin{document}
	
\maketitle

\begin{abstract}
	The theory of slow invariant manifolds (SIMs) is the foundation of various model-order reduction techniques for dissipative dynamical systems with multiple time-scales, e.g. in chemical kinetic models. The construction of SIMs and many approximation methods exploit the restrictive requirement of an explicit time-scale separation parameter. Most of those methods are also not formulated covariantly, i.e. in terms of tensorial constructions. We propose an intrinsically coordinate-free differential geometric approximation criterion approximating normally attracting invariant manifolds (NAIMs). We translate some ideas behind existing approximation approaches, the stretching based diagnostics (SBD) and the flow curvature method (FCM) to tensors of Riemannian geometry, specifically to spacetime curvature in extended phase space. For that purpose we derive from flow-generating smooth vector fields a metric tensor such that the original dynamical system is a geodesic flow on a Riemannian manifold. We apply the resulting method to test models.
\end{abstract}


\vspace*{4mm}
\AMS{37D99, 37M21, 53B50}

\section{Introduction}
A wide range of natural processes are modeled by high dimensional dynamical systems with multiple time-scales, for example in chemical kinetics. Their numerical treatment is challenging due to high dimension and stiffness resulting from spectral gaps. The existence of different time-scales usually correlates with a bundling behavior of solution trajectories near a lower-dimensional invariant manifold in phase space. By restriction to this manifold, both the curse of high-dimensionality and stiffness can be reduced significantly, resulting in suitable model-order reduction strategies.

The origins of invariant manifold theory reach back to the works of Lyapunov \cite{lyapunov1992general}, Hadamard \cite{hadamard1901iteration} and Poincar{\'e} \cite{poincare1899methodes}. Lyapunov's auxiliary theorem provides the existence and uniqueness of an analytic manifold tangent to the slow subspace in an equilibrium, as long as 'non-resonance' conditions are satisfied. The latter also guarantee the existence of invariant tori after non-linear perturbation of a system, according to the KAM-Theorem (see \cite{arnol1963proof}). Two popular concepts are normally hyperbolic invariant manifolds (NHIMs) as studied in \cite{hirsch2006invariant, wiggins2013normally,eldering2013normally}, and inertial manifolds (see \cite{temam2012infinite}). These notions are related: Inertial manifolds are normally hyperbolic in specific cases\cite{rosa1996inertial}. NHIMs also lay the mathematical foundation for Fenichel \cite{fenichel1971persistence, fenichel1974asymptotic,fenichel1977asymptotic, fenichel1979geometric} to use the Geometric Singular Perturbation Theory (GSPT) for slow invariant manifold (SIM) construction for singularly perturbed slow-fast systems.

 Quite a number of different approaches to compute low-dimensional manifolds for the purpose of model-order reduction have been developed. Some methods are directly rooted in chemistry \cite{bodenstein1913theorie,Chapman1913}. Others take the dynamical systems viewpoint such as the intrinsic low dimensional manifold (ILDM) \cite{maas1992simplifying}, the computational singular perturbation method (CSP) \cite{lam1989understanding,lam1994csp}, an iterative iterative method by Roussel and Fraser (RFM) \cite{roussel1991geometry,Roussel2012}, the G-scheme \cite{valorani2009g}, zero derivative principle (ZDP) \cite{gear2005projecting} and approaches \cite{lebiedz2010minimal,lebiedz2010entropy,lebiedz2011geometric,lebiedz2011variational,lebiedz2013continuation,lebiedz2014optimization,Lebiedz2016} using entropy and variational principles, just to name a few.

 Most of these approaches share underlying concepts. In \cite{Lebiedz2016}, e.g., it is shown that one class of methods utilizes derivatives of state vectors, another class a boundary value problem for trajectories. A number of methods approximate the SIM in GSPT for slow-fast system with different order of the asymptotic expansion. In particular ILDM with order one \cite{kaper2002asymptotic} as well as CSP \cite{zagaris2004analysis}, ZPD \cite{gear2005projecting} and RFM \cite{kaper2002asymptotic} with order proportional to the iteration/order of derivative. Both the CSP and ZDP generate suitable coordinate systems in the tangent bundle \cite{zagaris2005two} with a view on slow-fast decomposition.
 
  However, many approaches do not provide tensorial formulations and are restricted to slow-fast systems. If the defining quantities providing lower dimensional manifolds are tensors, one can choose coordinates at will which can be a significant benefit for application purposes. To the best of our knowledge, the only well-established method providing a tensorial formulation is the CSP as shown by Kaper et al. \cite{kaper2015geometry}.

 The central objectives of our work  are: (A) Find geometrically motivated tensorial formulations for established SIM methods either directly, or by translating them into a suitable setting. And (B): Reach goal (A) such that we can exploit the notion of normal hyperbolicity implying that we are not restricted to slow-fast systems. We approach these objectives by embedding dynamical systems in a setting similar to general relativity, interpreting their solutions as space-time geodesics in Riemannian geometry. This yields several advantages: The central quantities of this field are all covariant values, for example curvature-tensors. Those quantities yield geometric interpretations just like a couple of SIM methods provide their own geometrical meaning. If we successfully translate these methods to quantities in Riemann geometry we directly receive such a tensorial formulation. It turns out that ideas from the well-known stretching-based diagnostics (SBD) \cite{adrover2007structure,adrover2007stretching} by Adrover et.al. and the flow curvature method (FCM) \cite{ginoux2006differential, ginoux2008slow} can used as guidance, when a suitable framework is used. The SBD refers to normal hyperbolicity. Translating the SBD approach into such a tensorial setting is a major focus of this work.

 The paper is organized as follows: The foundations of normal hyperbolicity and GSPT as well as the benefit of tensorial formulations are discussed in \Cref{sec:NH}, \Cref{sec:GSPT} and \Cref{sec:adv_covariant} respectively. A suitable geometric setting is motivated and developed in \Cref{sec:Diffgeo_chapter}.
The resulting setting  provides the possibility to make use of various notions of \textit{intrinsic} curvature to analyze the bundling behavior of trajectories near the NAIMs. In \Cref{sec:geo_adrover}, we derive which of these notions are an evident choice by exploiting the geometrical foundation of the SBD \cite{adrover2007structure}. In \Cref{sec:geo_ginoux} we introduce new viewpoints on the FCM and illustrate, how it can reformulated in our Riemannian geometry setting. 

\subsection{Normal Hyperbolicity}\label{sec:NH}

 The concepts of normal hyperbolicity is central for the tensorial SIM approximation proposed in \Cref{sec:Diffgeo_chapter} and \Cref{sec:geo_adrover}. NHIMs have the defining property to admit a hyperbolic splitting into whitney sum of flow-invariant subbundles of the tangent bundle: 
	\begin{equation}\label{eq:whitney} T_M N = TM \oplus E^s \oplus E^u. \end{equation}
	Here, stable and unstable subbundles are indicated by $E^s$ and $E^u$ respectively. Since we are primarily interested in exponentially attracting manifolds. Hence, we restrict ourself to the case without unstable directions (though the following approach can be modified to suit the more general case). The sum \eqref{eq:whitney} reduces to $T_M N = TM \oplus E^s$ and $M$ is called normally hyperbolic attracting manifold (NAIM). Roughly speaking, normal hyperbolicity means that 
	\begin{enumerate}
		\item[(a)] The linearized flow along $M$ contracts along $E^s$ 
		\item[(b)] The contraction in $E^s$ dominates the dynamic in $TM$
	\end{enumerate}
	The latter property is formalized by considering so-called generalized Lyapunov-type numbers (see \cite{wiggins2013normally}), which is discussed in \Cref{sec:adrover_ori} in more detail.

\subsection{Geometric Singular Perturbation Theory }\label{sec:GSPT}

The established analytical foundation of the theory of SIMs is introduced in \cite{fenichel1971persistence,fenichel1974asymptotic,fenichel1977asymptotic,fenichel1979geometric} and can be applied to explicit slow-fast systems. This is a class of dynamical systems which can be written in the form
\begin{alignat*}{2}
\frac{\d}{\d t} x_s &= f(x_s,x_\f,\varepsilon)  \qquad &x_s &\in \R^{n_s} \\
\varepsilon\frac{\d}{\d t} x_{\f} &= g(x_s,x_\f,\varepsilon)\qquad &x_\f &\in \R^{n_\f}
\end{alignat*}
where $0 < \varepsilon \ll 1$. In the former setting, $x_s$ and $x_\f$ are called slow and fast variables respectively, such that $x=[x_s,x_\f]$ and $n_s+n_\f = n$. Fenichel showed (see \cite{fenichel1979geometric}) a persistence property under certain conditions: If there exists a NHIM for the vector field associated with $\varepsilon =0$, then there is also a NHIM for small $\varepsilon>0$. If a NHIM is constructed in this specific context in this particular way, we call it a SIM. In case of exponentially attracting invariant manifolds, a SIM is also a NAIM. In this case, a SIM is represented by the mapping 
\[
h_{\varepsilon}: \R^{n_s} \to \R^{n_\f}, \qquad x_s = h_{\varepsilon} (x_\f).
\]
where $h_{\varepsilon}$ can be expressed by a power series (asymptotic expansion) in $\varepsilon$:
\begin{equation} \label{eq:power-series}
h_\varepsilon (x_\f) = \sum_{k = 0}^{\infty} h_k(x_\f) \varepsilon^k 
\end{equation}
The functions $h_k(x_\f)$ in \cref{eq:power-series} are iteratively calculated by use of the so-called \textit{invariance equation} 
\begin{equation}\label{eq:invariance}
\varepsilon Dh(x_\f) f(x,h_{\varepsilon},\varepsilon) = g(x_\f,h_{\varepsilon}(x_\f),\varepsilon)
\end{equation} 
and matching of the coefficient with respect to $\varepsilon$-powers ({\it matched asymptotic expansion}).

Many of the computational methods referred to in the Introduction aim at an approximation of the $\varepsilon$-Taylor series to some order in some coordinate system.
Our aim is to find a local tensorial, coordinate independent, intrinsically geometric quantity approximating normal hyperbolicity itself, without restriction to singular perturbed slow-fast systems. We also obtain a SIM approximation for slow-fast systems as a mere consequence.

\subsection{Advantages of a covariant formulation}\label{sec:adv_covariant}
We use tensors fields as the fundamental concept in this work to propose a suitable covariant approximation method. Tensors are multilinear mappings from a cartesian product of vector spaces and dual spaces to $\R$. In our case, the vector space is the tangent space $T_p M$ of the solution manifold $M$ of the dynamical system $\dot{x}=f(x)$ in extended phase space (including a time coordinate). We base our approach on purely geometric concepts and choose our tensors accordingly: 
	\begin{itemize}
		\item[(a)] In \Cref{sec:ut_diffgeo} we introduce a specific metric (rank two tensor).
		\item[(b)] The metric is utilized to define the Riemann curvature tensor (rank four tensor) in \Cref{sec:deviation}.
		\item[(c)] By plugging-in tangent vectors (rank one tensors) we then receive the deviation tensor (rank two)  and a scalar value corresponding to a certain curvature (rank zero tensor) in \Cref{sec:geod_str}.
	\end{itemize}
For a given dynamical system $\dot{x} = f(x)$, the set of variables $X:=\{x_1, \ldots, x_n\}$ together with a time-variable $\tau$ forms a canonical coordinate-system for the space-time manifold $M$. All calculations in the above steps (a) - (c) are based on the bundling behavior of the solution trajectories of $\dot{x}=f(x)$ regarding the so-called 'parent coordinates' $X$ and a vector field $f$ expressed in these coordinates and thus defining the solution manifold geometrically.

The parent coordinates induce a basis $B_X$ of each tangent space (and cotangent space). Each tensor $\mathcal{T}$ can then be expressed by its coefficients in the $X$-coordinate frame. Conversely, once a basis and all coefficients regarding that basis are specified, $\mathcal{T}$ is well-defined as a coordinate-free object in the following manner:
Its action on every applicable combination of vectors and covectors does not depend on the choice of the basis.
 We can take another basis $B_Y$ (induced by a chart $Y$) and the coefficients of $\mathcal{T}$ transform in a certain manner, as do the coefficients of vectors and covectors we are inserting. This is what covariant means in our context: The definite way tensor coefficients change when transforming the basis or equivalently the coordinate system.

 The scalar curvature value in (see (c)) is used to determine the approximate NAIM location. Crucially, this curvature value is just the evaluation of a tensor and does not depend on the used basis. Hence, we can choose our own (so-called 'utilized') coordinate chart $Y$ by applying some transformation $y = \phi (x)$ and calculate steps (a) - (c) regarding $Y$. We receive the same tensors - still containing the same information of the original system $\dot{x}=f(x)$ - but now with different coefficients regarding $Y$. 
 
 This approach is applicable for any system of the form $\dot{x}=f(x)$ (as long as $f$ is smooth enough) and without having a presumed division into slow and fast variables. We can choose slow and fast variables or transform the coordinates entirely, according to our needs. 
 
 \textbf{Remark:} This approach requires us to originally choose 'parent coordinates' $X$ in which the dynamic system $\dot{x} = f(x)$ is expressed. All calculations - regardless of the choice of 'utilized coordinates' $Y$ - still contain the geometry regarding $X$. These computations are not invariant with respect to transformations, implying that our approximation is not invariant, too. This non-invariance is a necessity, since especially SIMs are not invariant to transformations as well. Consider the systems
 \[
 (\text{I}) \begin{cases}
 \dot{x}_1 &= -x_1 \\
 \dot{x}_2 &= -\varepsilon x_2
 \end{cases}
 \qquad \text{and} \qquad(\text{II}) \begin{cases}
 \dot{y}_1 &= -y_1 \\
 \dot{y}_2 &= -\frac{1}{\varepsilon} y_2
 \end{cases}
 \]
 with $0 < \varepsilon < 1$. System (I) can be transformed into (II) by 
\[
\begin{pmatrix}
y_1 \\ y_2
\end{pmatrix} = \phi(x_1,x_2) := \begin{pmatrix}
x_1\\
\left(\frac{1}{\varepsilon}\right)^2x_2
\end{pmatrix}
\] System (I) has a SIM at $\{x_1 =0\}$ and (II) at $\{y_2 = 0\}$ (not invariant regarding $\phi$). Our method - and basically  all other methods - can identify both SIMs separately. Covariance means that we can use $\phi$ as a local map and  $\{y_1, y_2\}$ as local coordinates in order to calculate the SIM of system (I). $\{x_1, x_2\}$ are 'parent coordinates' and $\{y_1, y_2\}$ are the utilized ones. In doing so, we receive a SIM at $y_2 =0$ which is incorrect with regard to system (II), but correct in terms of (I) by identifying $x_2 = \phi^{-1}(y) = 0 \Leftrightarrow y_2 =0$. Conversely, if we choose $\{y_1,y_2\}$ as parent coordinates, we receive different tensors and a SIM at $y_1 =0$ (correct for system (II)).
 
 In \Cref{sec:transf_rules}, we demonstrate how tensors and their coefficients change when choosing a different 'parent system' and how a covariant change of coefficients is calculated.

\section{Solution Trajectories as Geodesics in Spacetime}\label{sec:Diffgeo_chapter}
This section briefly introduces the main differential geometric setting of this work. It also discusses the motivation of choosing this specific setting based on geometric observations and physical analogies.
\subsection{Geometric Motivation}
In dissipative multiple time-scale systems solution initially fast trajectories for arbitrary initial values often converge towards an invariant submanifold while slowing down.

In extended phase space, by the introducing time $\tau$ as an additional axis, this behaviour is still observed.
\begin{definition} \label{def:ext_system}
	Let $f \in \mathcal{C}^{\infty}(E, \R^n)$, where $E \subset \R^n$ is an open set. We call the dynamical system \[
	\frac{\d}{\d t}
	\begin{pmatrix}
	x(t)\\
	\tau(t)
	\end{pmatrix}
	=
	\begin{pmatrix}
	\dot{x}\\
	\dot{\tau}
	\end{pmatrix}
	= \begin{pmatrix}
	f(x) \\
	1
	\end{pmatrix}, \qquad (x,\tau) \in E \times \R
	\]
	the extended system. 
\end{definition}

\begin{figure}[ht]
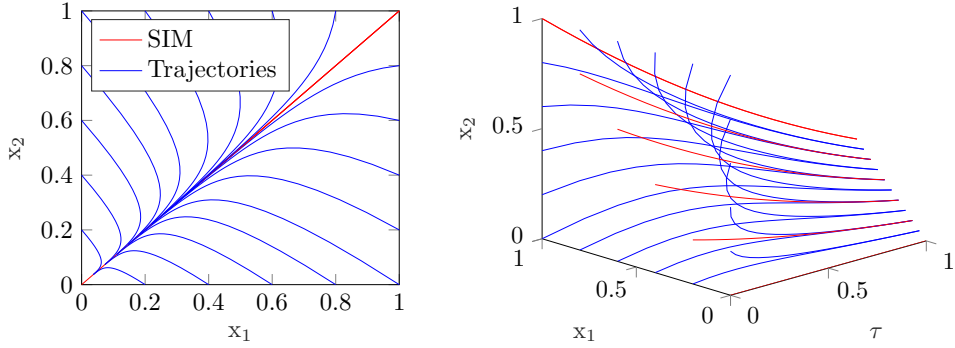
 \label{fig:bundling_plots}
\begin{minipage}[t]{.43\textwidth}
	\resizebox{1\textwidth}{!}{%
		\input{pictures/bundling_2d}
	}
\end{minipage}
\hspace{\fill}
\begin{minipage}[t]{.53\textwidth}
		\resizebox{1\textwidth}{!}{%
	\input{pictures/bundling_3d}
}
\end{minipage}
\caption{Phase space and extended phase space plot of linear system \cref{eq:lin_sys}. SIM in red.}
\end{figure}

\cref{fig:bundling_plots} illustrates how bundling behaviors of trajectories the original system relates to bundling of those of the extended one for the two-dimensional linear system
\begin{equation}\label{eq:lin_sys}
\frac{\d}{\d t}
\begin{pmatrix}
x_1 \\ x_2 
\end{pmatrix} = \left(
\begin{array}{cc}
-1- \gamma & \gamma \\ \gamma & -1-\gamma
\end{array} \right) \begin{pmatrix}
x_1 \\ x_2 
\end{pmatrix}
\end{equation}
for $\gamma = 3$ and different initial values. On the left plot, there are the solution trajectories of the original system, the SIM is curve (one-dimensional manifold). The right plot shows solution trajectories of the corresponding extended system. The SIM is a two-dimensional nonlinear surface spanned by solution trajectories embedded in $\R^3$.
\vspace*{2em}

The goal of this section is to formulate a setting capable of translating established SIM approximation methods into Riemannian geometry. Riemannian manifolds $(\M,g)$ - consisting of a smooth manifold $\M$ and a metric tensor field $g$ - are the foundations of that framework. We suggest a Riemannian manifold $(\M = \M_f,g=g_f)$ - which we call $f$-manifold - depending on a given vector field $f$ and show that the metric $g_f$ constructed from $f$ qualifies for an appropriate geometric setting.

 A recent work of Heiter and Lebiedz \cite{heiter2018towards} reformulates flow-invariance by vanishing specific time-sectional curvatures of submanifolds of the extended phase-space $\R^{n+1}$. Our work adopts from \cite{heiter2018towards} the idea of considering the extended phase-space $(x,\tau) \in \R^n \times \R$ in order to construct the metric $g_f$ with the desirable properties.
	
\subsection{Utilization of differential geometry}\label{sec:ut_diffgeo}
The field of Riemann geometry offers a wide variety of geometric quantities defined intrinsically, i.e. without reference to an embeeding space. We define a metric $g$ on the open set $E \times \R$ (which we call $\M$) giving rise to a Riemannian manifold $(\M,g)$. A detailed overview of Riemann geometry and curvature can be found e.g. in \cite{lee2006riemannian}. Our metric $g$ is chosen in a way that turns every solution trajectory into a geodesic - a shortest connection path with respect to the metric $g$. All differential geometric quantities used in this work depend on this specific metric $g$ which itself depends on the given dynamical system and is computed from the generating vector field. 

We integrate the former ideas into a mathematical formalism and introduce the  basic notions of differential geometry.
In the following definitions, we always assume that $f \in \mathcal{C}^{\infty}(E, \R^n)$ for some open set $E \subset \R^n$ and $n \in \mathbb{N}$ is fixed.
\begin{definition}
	The set $\mathcal{M}:= E \times \R$ defines a smooth manifold and the identity mapping
	\[\emph{id}: \M  \to \R^{n+1} \] is a local (and global) chart.
\end{definition}
We call the first $n$ coordinates of this chart $x_1(p),\ldots,x_n(p) = p_1, \ldots, p_n$ the state-components.
In contrast, the last coordinate $\tau(p) = p_{n+1}$ is the so-called time-component.

Let $T_p \mathcal{M}$ and $T_p^\prime \M$ denote the tangent space and cotangent space respectively for each point $p$. 
The set of derivatives in the direction of each coordinate forms a basis of $T_p\M$. These tangent vectors are denoted by
\[
\partial_{i,p} =\left. \frac{\partial}{\partial x^i}\right\vert_p, \qquad i = 1, ..., n \qquad \text{and}\qquad \partial_{n+1,p}:= \left. \frac{\partial}{\partial \tau}\right\vert_p
\]
The corresponding dual basis consisting of covectors is denoted by $\d x_{1,p}, ..., \d x_{n,p}$ and $\d \tau_p$. For $k,\ell\in \mathbb{N}$, $T_\ell^{k}(\M)$ indicates the set of all $k$-times covariant and $\ell$-times contravariant tensor fields on $\M$. $T_0^{1}(\M)= T \M$ and $T_1^0 \M$ represent the tangent bundle and cotangent bundle respectively.
We denote the base vector fields $\partial_{i}$ by 
\begin{align*}
\partial_i := \{ \partial_{i,p}~|~p \in \M \} \in T\M \qquad \forall i=1,...,n+1. 
\end{align*}
The base covector fields $\d x_i, \d \tau$ are defined in the same manner. 

We define the specific metric $g$ used in this work as a tensor field and state its basic properties: 
\begin{definition} \label{def:metric_tensor}
	Let $\M$ be as in \cref{def:ext_system}. Then, the mapping $g: \M \to T_2^0 (\M)$
	\begin{align*}
	p \mapsto g_p &= \left( \sum_{k = 1}^n  (\d x_{k,p} \otimes \d x_{k,p})  - f_k(x_p) ( \d \tau_p \otimes \d x_{k,p} + \d x_{k,p} \otimes \d \tau_p ) \right)\\
	&+ \left( 1+ \sum_{k = 1}^n f_k(x_p)^2  \right) \left( \d \tau_p \otimes \d \tau_p \right)
	\end{align*}
	defines a smooth tensor field $g$ on $\M$, where $\otimes$ indicates the tensor product. 
\end{definition}
For every fixed $p \in \M$, $g_p$ is represented by its components $g_{ij} = g_{p,ij} := g_p(\partial_{i,p},\partial_{j,p})$ with respect to the basis $\{\partial_{i,p}~|~i = 1,...,n+1\}$. For every $p \in \M$ and $v_p, w_p \in \T_p \M$ we have base representations
\[
v_p = \sum_{i=1}^{n+1} v^i \partial_{i,p}, \quad, \quad w_p =  \sum_{i=1}^{n+1} w^j \partial_{j,p}.
\]
The metric $g_p$ (as a symmetric bilinear form on tangent space) applied to the tuple $(v_p,w_p)$ then can be calculated by
\begin{align*}
g_p(v_p,w_p)&=
g_p\left( \sum_{i=1}^{n+1} v^i \partial_{i,p}, \sum_{i=1}^{n+1} v^i \partial_{j,p} \right)\\ &= \sum_{i,j = 1}^{n+1} v^i w^j g_p(\partial_{i,p},\partial_{j,p})
= v^T \left(g_{ij}\right)_{1\le i,j \le n+1} w.
\end{align*}
The components of $g_p$ can be deduced from \cref{def:metric_tensor} and read
\begin{equation} \label{eq:metr_comps}
\left(g_{p,ij}\right)_{1\le i,j \le n+1} = \begin{pmatrix}
\text{Id}_n & -f(x_p)\\
-f(x_p)^T & 1+\| f(x_p) \|_2^2
\end{pmatrix},
\end{equation}
where $\text{Id}_n$ indicates the $n \times n$ identity matrix.
\begin{proposition}
	Let $\M$ be defined as before and $p = [x_p,\tau_p] \in \mathcal{M} $ an arbitrary point. The tensor $g_p$ is a metric for every fixed $p$, independent of the values of $f(x_p) \in \R^n$.
\end{proposition}

\begin{proof}
	The tensor $g_p$ is a bilinear form at every point $p$ by the definition of the tensor product. The symmetry of the matrix $\left(g_{p,ij}\right)_{1\le i,j\le n+1}$ in equation \cref{eq:metr_comps} implies the pointwise symmetry of $g_p$. 	
	It suffices to show that the matrix $\left(g_{p,ij}\right)_{1\le i,j\le n+1}$ is positive definite for every value of $f(x(p))$. Since the identity $\text{Id}_n$ is positive definite implying that all its minors are positive, we only have to show that $\det (g_{p,ij})>0$. Adding $f_i(x_p)$-times the $i$-th column to the $(n+1)$st column for each $1\le i \le n$ yields
	\[
	\det \begin{pmatrix}
	\text{Id}_n & -f(x_p) \\
	-f(x_p)^T & 1+\|f(x_p)\|_2^2
	\end{pmatrix}  = \det \begin{pmatrix}
	\text{Id}_n & 0 \\
	-f(x_p)^T & 1
	\end{pmatrix} = 1.
	\]
	Hence, the matrix $g_{p,ij}$ is positive definite and $g_p$ is a metric tensor for each $p \in \M$. 
\end{proof}

\begin{corollary} \label{cor:Rie_mani}
	For any given smooth function $f$, the tuple $(\M,g)$ is a Riemannian manifold.
\end{corollary}
We call the tuple $(M,g)$ f-manifold.
The right hand side of the extended system in \cref{def:ext_system} defines a smooth vector field $\T:\M \to T\M $ on $\M$. Its coordinate representation is given by 
\begin{equation} \label{eq:tangent_vec}
\T_p =\T(p) = \sum_{k = 1}^n f_k(x_p) \partial_{k,p} + \partial_{n+1,p} \qquad \forall p \in \M
\end{equation}
The extended system is a dynamic system on $\M$. The core property of the metric $g$ is formalized in the following theorem:
\begin{theorem} \label{thm:geodetization}
	Let $f:E \to \R^n$ be given, $(\M,g)$ be the corresponding f-manifold and $\nabla =\nabla^g$ be the Levi-Civita connection which preserves $g$. Let $\gamma: (-\varepsilon, \varepsilon) \to \M$ be a solution curve of the extended system of $f$ on $\M$. Then, $\gamma$ is a geodesic with respect to $\nabla$. In particular, $\gamma$ satisfies the geodesic equation 
	\[
	\nabla_{\dot{\gamma}}\dot{\gamma} = 0
	\]
	at every point $\gamma(t) = [x(t),\tau(t)] \in \M$. With regard to the coordinates $(x,\tau)$, the former equality reads
	\begin{equation}\label{eq:geod_eq}
	\frac{\d^2}{\d t^2}\begin{pmatrix}
	x(t) \\ \tau(t)
	\end{pmatrix} =- \left(  \left(\frac{\d}{\d t}( x(t), \tau(t)) \right) (\Gamma_{ij}^k(\gamma(t)))_{i,j} \frac{\d}{\d t} 
	\begin{pmatrix}
	x(t)\\
	\tau(t)
	\end{pmatrix} \right)_{k = 1,...,n+1},
	\end{equation}
	with $ \Gamma_{ij}^k(\gamma(t)) $ being the Christoffel symbols of the Levi-Civita connection $\nabla$ evaluated at the point $\gamma(t)$.
\end{theorem}

\begin{proof}
	By calculation of the Christoffel symbols $\Gamma_{ij}^k$, see Appendix.
\end{proof}
Every solution trajectory of the extended system has equal velocity, since
\[
g_p(\T_p,\T_p) =  [f(x_p)^{\text{T}},1] \begin{pmatrix}
\text{Id}_n & -f(x_p)\\
-f(x_p)^T & 1+\| f(x_p) \|_2^2
\end{pmatrix} \begin{pmatrix}
f(x_p)\\1
\end{pmatrix} = 1.
\]
In this sense, the metric is a normalizer of the time parametrization of solution trajectories. 

The proposed setting shares similarities with general relativity where the trajectories of free-falling particles are geodesics with regard to a metric representing a gravitational field. This interpretation motivates the use of concepts from general relativity - such as geodesic deviation - to represent attractiveness. This is new in the SIM context and the idea is implemented in \Cref{sec:geo_adrover}.
\subsection{Covariant Transformation of the metric} \label{sec:transf_rules}
We briefly demonstrate how the coefficients of the metric tensor $g$ from definition \ref{def:metric_tensor} change under simple re-scaling of parent coordinates. Let $\dot{x}=f(x)$ be the parent system with coordinates $X:=\{x_1,\ldots,x_n, \tau\}$. The metric tensor from the last section regarding this choice of parent coordinates is denoted by $ \left[ g \right]_{(1)}$. Its coefficients regarding $X$ are denoted by $\left[ g_{ij}^{(X)} \right]_{(1)}$ (see \ref{eq:metr_comps}). Consider new coordinates $Y := \{y_1, \ldots, y_{n+1}\}$ (where $y_{n+1}$ becomes the new explicit time coordinate) obtained by the transformation
\begin{equation}\label{eq:simple_transf}
\begin{cases}
y_k = \phi_k(x,\tau) :=  a_k x_k \qquad \text{with} \qquad a_k \ne 0 \quad\forall k = 1, \ldots,n.\\
y_{n+1} = \phi_{n+1}(x,\tau) := a_{n+1} \tau \qquad \text{with} \qquad a_{n+1} = 1 
\end{cases}
\end{equation}
In $Y-$coordinates, the phase-space system then reads $\dot{y}_k = (a_k f(x(y)))_{k=1,\ldots,n}$. We can express $\left[ g \right]_{(1)}$ by the means of $Y$ and calculate
\begin{align*}
\left[ g_{ij}^{(Y)} \right]_{(1)} :&= \left[ g \right]_{(1)}\left(\frac{\partial}{\partial y_i},\frac{\partial}{\partial y_j}\right) = 
\left[ g \right]_{(1)}\left( \sum_{k=1}^{n+1} \frac{\partial \phi_k^{-1}} {\partial y_i} \frac{\partial}{\partial x_k} , \sum_{\ell=1}^{n+1} \frac{\partial \phi_\ell^{-1}}{\partial y_j} \frac{\partial}{\partial x_\ell} \right)\\
&= \sum_{k,\ell=1}^{n+1} \frac{\partial \phi_k^{-1}} {\partial y_i} \frac{\partial \phi_\ell^{-1}}{\partial y_j}  \left[ g_{ij}^{(X)} \right]_{(1)} 
= \sum_{k,\ell=1}^{n+1} \delta_{ki} \frac{1}{a_i} \delta_{\ell j} \frac{1}{a_j}  \left[ g_{ij}^{(X)} \right]_{(1)} \\
&= \frac{1}{a_i a_j} g_{ij}^{(X)} \qquad \forall (i,j) \in \{1, \ldots, n+1\}^2
\end{align*}
Inserting the transformation \ref{eq:simple_transf} and original coefficients from definition \ref{def:metric_tensor}, we receive
	\[
\left[g_{ij}^{(Y)}\right]_1 = \left(
\begin{array}{c c}
A_{-2} & (-\frac{1}{a_k} f_k)_k\\  (-\frac{1}{a_k} f_k)_k^T &1 + \|( f_k)_k\|_2^2\\
\end{array}
\right)
\qquad \text{with}  \qquad A_{-2} := \begin{pmatrix}
\frac{1}{a_1^2} & \ldots & 0 \\
\vdots & \ddots& \vdots \\
0 & \ldots & \frac{1}{a_n^2}
\end{pmatrix}.
\]
Alternatively, when declaring $\dot{y}_k = (a_k f(x(y)))_k$ as parent system, we receive a different metric $ \left[ g \right]_{(2)}$. The coefficients have to be chosen according to definition \ref{def:metric_tensor}. Hence, the resulting coefficients with respect to $Y$ are
\[
\left[g_{ij}^{(Y)}\right]_2 =\left(  \begin{array}{c c}
	\text{Id}_n &(-a_k f_k)_k \\
	(-a_k f_k)_k^T & 1 + \|(a_k f_k)_k\|_2^2
\end{array}
\right) \ne \left[g_{ij}^{(Y)}\right]_1
\]
The metrics are different, thus all derived tensors as well as our calculated manifold. 

\section{Geodesic Stretching Approach} \label{sec:geo_adrover}
In this section, we introduce a new tensorial method to approximate normal attractiveness by exploiting the previously introduced setting in the following way: We translate the notion of geodesic deviation to the concept of the ration of stretching rates from SBD (see \cite{adrover2007structure,adrover2007stretching}). The results are what we call geodesic stretching rates which turn out to be specific sectional curvatures, in coordinates corresponding to some curvature tensor entries. In \ref{sec:test_stret} we apply the resulting method to non-linear test-models. In \ref{sec:test_stret} we derive a SIM approximation method from this approach and apply the resulting method to non-linear test-models.

\subsection{Deviation}\label{sec:deviation}
In general relativity geodesic deviation is used to describe relative behavior of neighboring geodesics corresponding to the relative acceleration of nearby particles in free-fall. It is defined by plugging in a tangent vector $y_p$ - representing the instantaneous velocity of the geodesic - into the first and third argument of the Riemann curvature tensor which is denoted by 
\[
\mathcal{R}_p: \left(T_p \M \right)^3 \to T_p \M \qquad (u_p, v_p, w_p) \mapsto \mathcal{R}(u_p, v_p)w_p \in T_p \M  \qquad\forall p \in \M.
\]
The result is a tensor field depending on the tangent vectors $y_p$. The input vector $v_p$ of this reduced tensor represents a small displacement between the neighboring geodesics, while the output stands for the deviation.
On the f-manifolds $(\M,g)$ from the previous \Cref{sec:Diffgeo_chapter}, there is one set of geodesics of special interest: The solution trajectories of the extended system bundling near a NAIM. Hence, an evident choice is $y_p = \T_p $ for all $p \in \M$ and receive tensor-field  depending on $\T_p$, leading to the following definition:
\begin{definition} \label{def:system_dev}
	Let $(\M, g)$ be as in \cref{cor:Rie_mani} and $\T$ as in equation \cref{eq:tangent_vec}. Let $\mathcal{R} = \mathcal{R}_f$ be the corresponding Riemann curvature tensor. We call the tensor field
	\[
	S \in T_1^1(\M), \qquad T_p \M \ni v_p \mapsto S(v_p) := R_p(\T_p, v_p) \T_p \in T_p \M\qquad p \in T_p \M
	\]
	f-deviation.
\end{definition}

\begin{remark}
	The christoffel symbols, the curvature tensor and the f- deviation do not depend on the explicit time $\tau$.
\end{remark}
\begin{proof}
	The components of $g$ are independent of the time $\tau$, implying that the components $g^{ij}$ of the inverse metric tensor are also time-independent. Since we use the Levi-Civita connection, the Christoffel symbols are calculated by derivatives of time-independent quantities and the statement holds for $\Gamma_{ij}^k$. Using the same argumentation, we conclude this property for the curvature tensor $\mathcal{R}$ and the f-deviation $S$.
\end{proof}
Based on its properties and geometric interpretations, the f-deviation appears to be well-suited to be turned into a geometric criterion to approximate normal attractiveness. We now aim to deduce a scalar, curvature-based quantity from the f-deviation that intuitively represents the bundling behavior.
In order to do so, we are guided by an existing, geometric approach to characterize normal attractiveness: The SBD, introduced in \cite{adrover2007structure}, \cite{adrover2007stretching} by Adrover et al. .

\subsection{Original Stretching approach}\label{sec:adrover_ori}
  SBD in dissipative and chaotic system is a local, geometric reduction approach to multiple time scale dynamics. Let $M$ be an embedded submanifold of $\R^n$ with $x$ being the euclidean coordinates and $p \in M$ an arbitrary point. Normal attractiveness requires vector dynamics in tangential direction to be dominated by the dynamic in normal direction. This property is encoded in the generalized Lyapunov-type number (see \cite{kuehn2015multiple} and \cite{wiggins1994normally})
  \begin{equation}\label{eq:GLTN}
  \sigma(p) = \inf \left\{ b~:~\frac{\|v_{-t}\|}{\|w_{-t}\|^b}\to 0\text{, as} ~t\to \infty,~\forall v_0 \in T_pM,~w_0\in N_PM\right\}
  \end{equation}
Here, $v_t = v(x(t))$ and $w_t = w(x(t))$ represent the vector dynamic along a solution of the ODE $\dot{x} = f(x)$ for tangential and normal initial vectors respectively. Combined with demanding $\|v_{-t}\|$ to approach infinity at a certain rate, normal attractiveness requires $\sigma(p)<1$. Hence, we compare the convergence rates of the limits 
\[
\lim_{t \to \infty} \|v_{-t}\|\quad \text{and} \quad \lim_{t \to \infty} \|w_{-t}\|.
\]
According to \cite{adrover2007stretching}, an evident local quantity to approximate contraction rates of vector dynamics in time, is to take the first derivative of the square: 
\[
\frac{\d \|v\|^2}{\d t} = 2 \langle J_f(p)v , v\rangle
\]
where the brackets $\langle \cdot, \cdot \rangle$ represent the euclidean inner product.
Normalizing the previous term by division by $2\|v\|^2$ - implying that the resulting terms do not depend on the vector's length - leads to the definition of so-called stretching rates $\omega_p(T)$ and $\omega_p(N)$ for tangential and normal vectors respectively, where
\begin{align*}
\omega_p(T) := \max_{w \in T_pM}&\frac{\langle J_f(x)  v, v\rangle}{\langle v ,v\rangle}, \\
\omega_p(N) :=\max_{w \in N_pM} &\frac{\langle J_f(x)  w, w\rangle}{\langle w ,w\rangle}.
\end{align*}
Instead of considering the limit of the quotient $(\|w_{-t}\|/\|v_{-t}\|)$, Adrover et. al. take the quotient of the stretching rates instead. These are so-called stretching ratios 
$
r(p) = \omega_p(N)/\omega_p(T)
$
and replace the generalized Lyapunov-type numbers $\sigma(p)$. We receive local geometric quantities representing NAIM properties: Normal stretching $\omega_p(N)$ represents bundling behavior near $M$, tangential stretching $\omega_p(T)$ indicates acceleration/deceleration alongside $M$, while their ratio $r(p)$ locally approximates normal attractiveness. An evident way of introducing a lower-dimensional manifold for model-order reduction is to maximize/minimize the stretchting rates and/or the corresponding ratio.

\subsection{Geodesic stretching}\label{sec:geod_str}
The stretching rates incorporate a geometric interpretation which is adopted to be transferred in the differential geometric setting from \Cref{sec:Diffgeo_chapter}. 
We can interpret $(\R^n, \langle \cdot, \cdot \rangle)$ as a Riemann manifold, equipped with euclidean metric. 
The term $J_f v$ represents the so-called vector dynamics of the flow $\phi(t)$ related to the differential equation $\dot{x} = f(x)$. The flow differential $D\phi_t: T_{x(0)} \R^n \to T_{x(t)} \R^n$ propagates perturbation vectors $v_{x(0)}$ along $\phi_t$. Hence the map $D\phi_t$ indicates how solution trajectories of $\dot{x}=f(x)$ diverge or converge. Aiming to extract a local rate of deviation, a direct calculation yields
\[
\lim_{t\to 0} \left( \frac{D\phi_t(v) -v}{t}  \right) = J_f v.
\]
Hence, the mapping $ J: T_x \mathbb{R}^n \to T_x \R^n, v \mapsto  J_f v $ assigns a perturbation vector $v$ to a local rate of deviation. This interpretation shares major similarities with the f-deviation introduced in \cref{def:system_dev}. 
An evident transfer of the $\omega_x(v)$ into the coordinate-free setting from \Cref{sec:Diffgeo_chapter} is to replace $(\R^n,\langle \cdot, \cdot, \rangle)$ by $(\M,g)$ and apply the f-deviation instead of $J$. The result is the so-called geodesic stretching rate:  
\begin{definition} \label{def:geod_stretching}
	Let $(\M,g)$ be defined as in \Cref{sec:Diffgeo_chapter}, $S$ the f-deviation, $p \in \M$ an arbitrary point and $v_p \in T_p \M$.
	The mapping 
	\[
	\vartheta_p: T_p \M \to \R, \qquad v_p \mapsto \frac{g_p(S_p(v_p),v_p)}{g_p(v_p,v_p)} \qquad \forall v_p \in T_p \M
	\]
	is called geodesic stretching. The image $\vartheta_p(v_p)$ is denoted as geodesic stretching rate of $v_p$.
\end{definition}
\begin{figure}[ht] \label{fig:vis_geo_stret}
	\centering
	\def\svgwidth{0.8\textwidth}
	\fbox{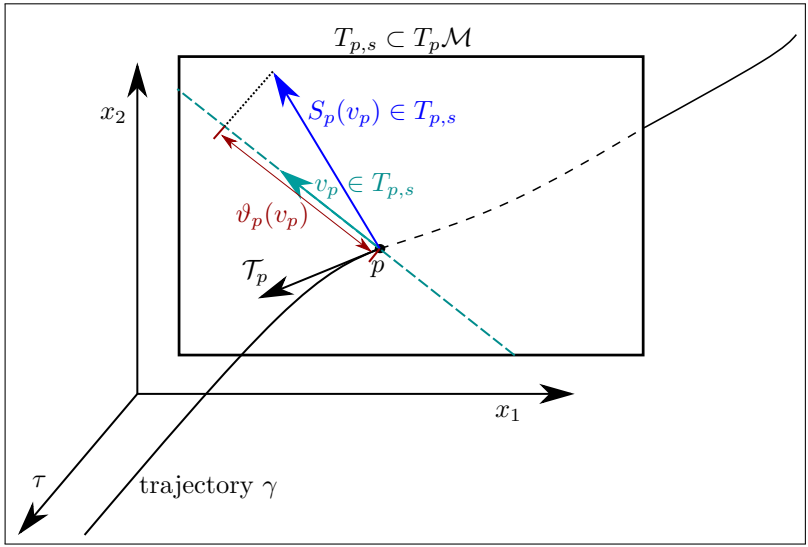}
	\caption{Visualisation of geodesic stretching rates}
\end{figure}
\begin{remark}
	The quantity $\vartheta_p(v_p)$ does not depend on the length of $v_p$ and is independent of explicit time $\tau$.
\end{remark}
By definition, we can calculate geodesic stretching rates for every tangent vector $v_p \in T_p \M$, especially those with non-vanishing time component $d\tau_p(v_p) \ne 0$. We exclude the latter ones and only consider those of the following subspace:
\begin{definition} \label{def:subspace}
	Let $p \in \M$ be arbitrary. We call the subspace
	\[
	T_{p,s} :=\{v_p \in T_p \M~|~\d \tau_p(v_p) = 0 \} = \text{span}(\partial_{1,p}, \ldots, \partial_{n,p})
	\]
	the pure-state space of $p$.
\end{definition}
This restriction appears natural, since the time axis $\tau$ is included artificially and we have no valid interpretation for "perturbation in explicit time". A visualization of the geodesic stretching rates can be found in \cref{fig:vis_geo_stret}.

There exists a curvature-based correspondent $\vartheta_p(v_p)$ that is well-defined for $v_p \in T_{p,s}$, formalized in the following theorem:
\begin{theorem}
	Let $(M,g)$, the stretching rates $\vartheta_p$, the subspace $T_{p,s}$ defined as above, $v_p \in T_{p,s} $ arbitrary. Then, $\vartheta_p(v_p)$ equals the sectional curvature of the subspace $\sigma_p^v$ spanned by the vectors $\T_p$ and $v_p$.
\end{theorem}
\begin{proof}
	Let $p \in \M$ and $v_p \in T_{p,s}$ be arbitrary. The set $\{\partial_{1,p}, \ldots \partial_{n,p}, \T_p \}$ is an orthonormal basis of $T_p \M$ for all $p\in \M$ with respect to $g_p$, implying $g_p(\T_p,v_p) = 0$. Using the fact that $g_p(\T_p,\T_p) = 1 $ we can calculate
	\[
	\vartheta_p(v_p) = \frac{g_p(R(v_p),v_p)}{g_p(v_p,v_p)} = \frac{g_p(\mathcal{R}_p(\T_p,v_p)\T_p,v_p)}{g_p(v_p,v_p)g_p(\T_p,\T_p)- g_p(\T_p,v_p)^2} =\mathcal{K}_p(\sigma_p^{v}),
	\]
	where $\sigma_p^{v} := \text{span}(\T_p,v_p)$ and $\mathcal{K}_p(\sigma_p^{v})$ is the sectional curvature of $\sigma_p^v$.
\end{proof}
\begin{corollary}
	The geodesic stretching rate $\vartheta(v_p)$ is a covariant intrinsically geometric quantity for every $v_p \in T_{p,s}$ and $p \in \M$. 
\end{corollary}

The aim is to exploit the rates $\vartheta(v_p)$ in an analogous way as the stretching characterize NAIMs by decomposing the tangent space into tangential and normal directions of a submanifold $U \subset \R^n$. Since $\M$ is a space-time manifold, we select the pure-state space $T_{p,x}$ and split it in the same manner as in the previous subsection. The result is the following definition:
\begin{definition} \label{def:pure_space_division}
	Let $x\in \R^n$ and $U \subset \R^n$  be an embedded submanifold with $x\in U$. Assume $\dim(T_x U) = k$ and vectors $u_{1,x}, \ldots u_{n,x} \in T_x \R^n $ satisfying \[ T_x U = \text{span} (u_{1,x}, \ldots, u_{k,x}),\qquad  (T_x U)^{\bot} = \text{span}(u_{k+1,x}
	\ldots,u_{n,x}).\]
	Let $[u_{1j}, \ldots , u_{nj}] $ be the euclidean coordinates of $u_{x,j}$ for all $j \in \{1,\ldots,n\}$, $p = [x,\tau] \in \M$ for an arbitrary $\tau$ and $x$ from above. 
	We then define the projected tangent space $T_{p,s}^{\text{tan}}$ and projected normal space $T_{p,s}^{\text{orth}}$ by
	\begin{align*}
	T_{p,s}^{\text{tan}}(U) &:= \left\{ \sum_{j=1}^{n} u_{ij} \partial_{j,p} ~|~ i = 1, \ldots,k \right\} \\ T_{p,s}^{\text{orth}}(U) &:= \left\{ \sum_{j=1}^{n} u_{ij} \partial_{j,p} ~|~ i = k+1, \ldots,n \right\}
	\end{align*}	
\end{definition}
By definition, we get $T_{p,s} = T_{p,x}^{\text{tan}}(U) \oplus T_{p,s}^{\text{orth}}(U)$.
\begin{definition} \label{def: tangVSorth_stret}
	Let $U$ be a submanifold of $\R^n$, $T_{p,s}^{\text{tan}}(U)$ and $T_{p,s}^{\text{orth}}(U)$ defined as above. 
	We define the tangential and orthogonal stretching rate as
	\[
	\varTheta_p^{\text{tan}}(U) := \max_{v_p \in T_{p,x}^{\text{tan}}(U)} \left( \vartheta_p(v_p)\right) \qquad \text{and} \qquad  \varTheta_p^{\text{orth}}(U) := \max_{v_p \in T_{p,x}^{\text{orth}}(U)} \left( \vartheta_p(v_p)\right)
	\]
	respectively.
\end{definition}
 Geodesic stretching rates are now used similarly to the original stretching rates in \cite{adrover2007stretching}: Tangential- and normal geodesic stretching rates measure deceleration and contraction respectively, while their ratio approximates the property of normal attractiveness.

A consequence of the geodesic stretching rate construction is the fact that it is capable of providing an approximation for a specifically constructed NAIM: the SIM in the context of slow-fast systems which is subject of the following section. We call the resulting method Geodesic Stretching method (GSM).
\subsection{Testing the Geodesic Stretching method} \label{sec:test_stret}
We conclude this section by applying the above method non-linear test models. We go into more detail on how this construction is deployed to approximate one-dimensional SIMs in these two-dimensional models. Finally, we discuss how the approach can then be applied to more general cases.
\subsubsection*{Geodesic Stretching for the Davis-Skodje Test Model}
 \begin{figure}[htp]
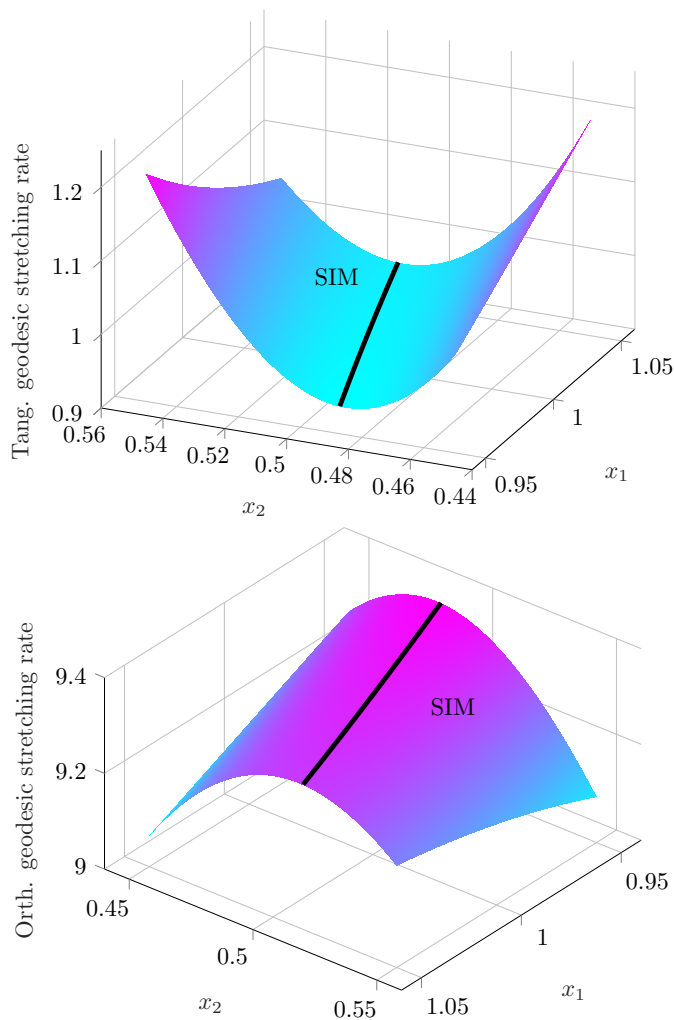
 \label{fig:tang_stret}
 	\centering
 	\resizebox{0.7\textwidth}{!}{%
 		
 		\input{pictures/stretching_tang_matlab}
 	}
 	\vspace*{\floatsep}
 	\resizebox{0.7\textwidth}{!}{
 		\input{pictures/stretching_orth_matlab}
 	}
 	\caption{Tangential geodesic stretching (first plot) and Orthogonal geodesic stretching (second plot) for the Davis-Skodje model; stretching on the SIM in black.}
 \end{figure}
Consider the non-linear Davis-Skodje system (see \cite{davis1999geometric})
\begin{subequations}
\begin{alignat}{2}
\dot{x}_1 &= -x_1 &=: f_1(x_1,x_2) \\
\dot{x}_2 &= -\eta x_2 + \frac{(\eta - 1)x_1+ \eta x_1^2}{(1+x_1)^2} &=: f_2(x_1,x_2)
\end{alignat}
\end{subequations}
where the parameter $\eta>1$ measures time-scale separation. 
 This system has a one-dimensional SIM with graph representation
\[
x_2 = h(x_1) = \frac{x_1}{1+x_1} \qquad \forall x_1 \in \mathbb{R}^+.
\]
For a two-dimensional system, the only SIM candidates are one-dimensional submanifolds $U$, i.e. the solution curves of the original system $\dot{x}= f(x)$.

\begin{figure}[htb] \label{fig:slice_stret} 
	\centering
\begin{tikzpicture}

\begin{groupplot}[group style={group size=1 by 2}, height = 0.4\textwidth, width=0.8 \textwidth]
\nextgroupplot[
ylabel = geodesic stretching,
legend cell align={left},
legend style={at={(0.03,0.97)}, anchor=north west, draw=white!80.0!black},
tick align=outside,
tick pos=left,
x grid style={white!69.01960784313725!black},
xmin=0.4945, xmax=0.5055,
xticklabel style={/pgf/number format/zerofill, /pgf/number format/precision=3},
xtick style={color=black},
y grid style={white!69.01960784313725!black},
ymin=0.947268662176504, ymax=0.950556716590399,
yticklabel style={/pgf/number format/zerofill, /pgf/number format/precision=3},
ytick style={color=black}
]
\addplot [semithick, cyan, mark=triangle*, mark size=3, mark options={solid}]
table {%
0.495 0.948147536662798
0.4955 0.947941761740631
0.496 0.947769956106313
0.4965 0.947632053355279
0.497 0.947527986397072
0.4975 0.947457687460525
0.498 0.947421088098941
0.4985 0.947418119195317
0.499 0.947448710967581
0.4995 0.947512792973833
0.5 0.947610294117647
0.5005 0.947741142653355
0.501 0.947905266191358
0.5015 0.948102591703469
0.502 0.94833304552826
0.5025 0.948596553376435
0.503 0.948893040336208
0.5035 0.949222430878713
0.504 0.949584648863419
0.5045 0.949979617543561
0.505 0.950407259571585
};
\addlegendentry{tang. geod. stret.}
\addplot [semithick, black]
table {%
0.5 0.947268662176504
0.5 0.950556716590399
};
\addlegendentry{SIM pos.}

\nextgroupplot[
ylabel = geodesic stretching,
xlabel = $x_2$ value,
legend cell align={left},
legend style={at={(0.03,0.03)}, anchor=south west, draw=white!80.0!black},
tick align=outside,
tick pos=left,
x grid style={white!69.01960784313725!black},
xmin=0.4945, xmax=0.5055,
xticklabel style={/pgf/number format/zerofill, /pgf/number format/precision=3},
xtick style={color=black},
y grid style={white!69.01960784313725!black},
ymin=9.3306932834096, ymax=9.3339813378235,
yticklabel style={/pgf/number format/zerofill, /pgf/number format/precision=3},
ytick style={color=black}
]
\addplot [semithick, purple, mark=asterisk, mark size=3, mark options={solid}]
table {%
0.495 9.3331024633372
0.4955 9.33330823825936
0.496 9.33348004389369
0.4965 9.33361794664473
0.497 9.33372201360293
0.4975 9.33379231253948
0.498 9.33382891190105
0.4985 9.33383188080469
0.499 9.33380128903242
0.4995 9.33373720702617
0.5 9.33363970588235
0.5005 9.33350885734664
0.501 9.33334473380865
0.5015 9.33314740829653
0.502 9.33291695447174
0.5025 9.33265344662357
0.503 9.33235695966379
0.5035 9.33202756912128
0.504 9.33166535113658
0.5045 9.33127038245644
0.505 9.33084274042841
};
\addlegendentry{orth. geod. stret.}
\addplot [semithick, black]
table {%
0.5 9.3306932834096
0.5 9.3339813378235
};
\addlegendentry{SIM pos.}
\end{groupplot}

\end{tikzpicture}
	\caption{Tangential and Orthogonal geodesic stretching rates for the Davis-Skodje model with $\eta = 3$ for $x_1 = 1$.}
\end{figure}
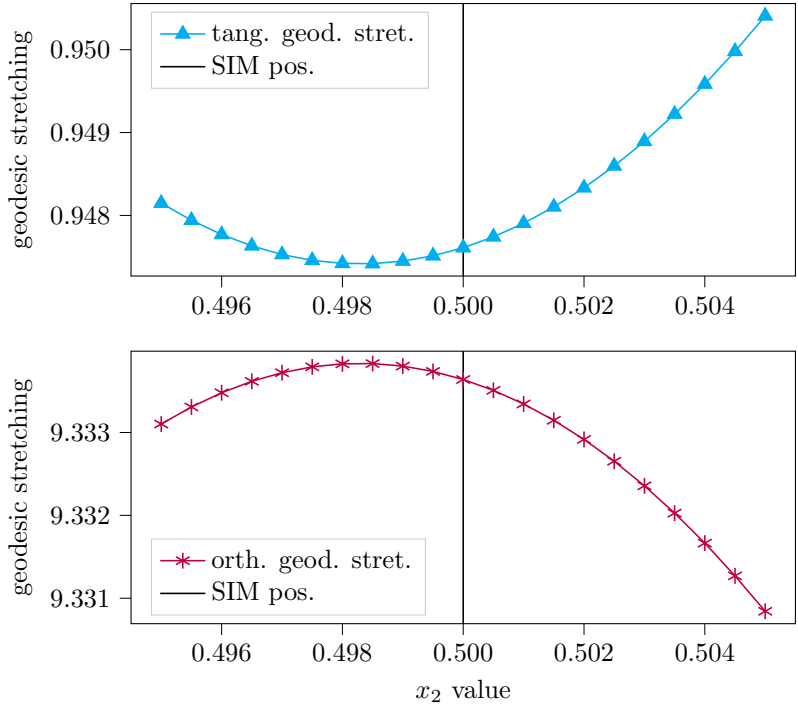
Hence, both the tangential and normal space of each trajectory is one-dimensional as well. Following \cref{def:pure_space_division}, we get the subspaces
\begin{alignat*}{2}
T_{p,s}^{\text{tan}}(\gamma) & = \text{span}(v_{1,p} ) ,\qquad v_{1,p}&:=f_1(x_p) \partial_{1,p}+f_2(x_p) \partial_{2,p}\\
T_{p,s}^{\text{orth}}(\gamma)  &= \text{span}(v_{2,p} ) ,\qquad v_{2,p}&:=f_2(x_p) \partial_{1,p}-f_1(x_p) \partial_{2,p}
\end{alignat*}
for each point $p$ in space-time $\M$. Since the subspaces are one-dimensional and the geodesic stretching rate does not depend on the length of each vector, we get
\[ \varTheta_p^{\text{tan}}(\gamma) = \vartheta_p(v_{1,p})\qquad \text{and} \qquad
\varTheta_p^{\text{orth}}(\gamma) = \vartheta_p(v_{2,p}).
\]
\cref{fig:tang_stret} depicts tangential and respectively orthogonal stretching rates in the vicinity of the SIM. 
Very close to the SIM, the tangential stretching rate gets small, while the orthogonal one gets particularly large. This observation matches the description of the SBD, introduced in the beginning of this section. By fixing one variable (e.g. $x_1$) and maximizing/minimizing orthogonal/tangential geodesic stretching rates with respect to the other variable, we get at least an adequate approximation of the SIM. 

Figure \ref{fig:slice_stret} shows that the former criterion is not exact. There, we fix $x_1 = 1$ and consider both tangential and orthogonal geodesic stretching as a function of $x_2$. Both resulting one-dimensional graphs have an extremum at around $x_2 = 0.4985$, while the SIM point is at $0.5$. We directly conclude that the ratio between both rates is also extremal at around $x_2 = 0.4985$.

\subsubsection*{ Comparing the Geodesic Stretching Method}
In this subsection, we refer to the chemical reaction mechanism
\begin{equation}\label{eq:c_MM_chemical}
\begin{cases}
\tilde{A}_1+ \tilde{A}_2 \rightleftharpoons \tilde{A}_3\\
\tilde{A}_3 \rightleftharpoons \tilde{A}_2+ \tilde{A}_4
\end{cases}
\end{equation}
taken from \cite{chiavazzo2007comparison} aiming to compare our method to different well-estabished methods to approximate the SIM. In this model, we have four species with respective concentrations $c_1, c_2, c_3$ and $c_4$. Apart from the ILDM,
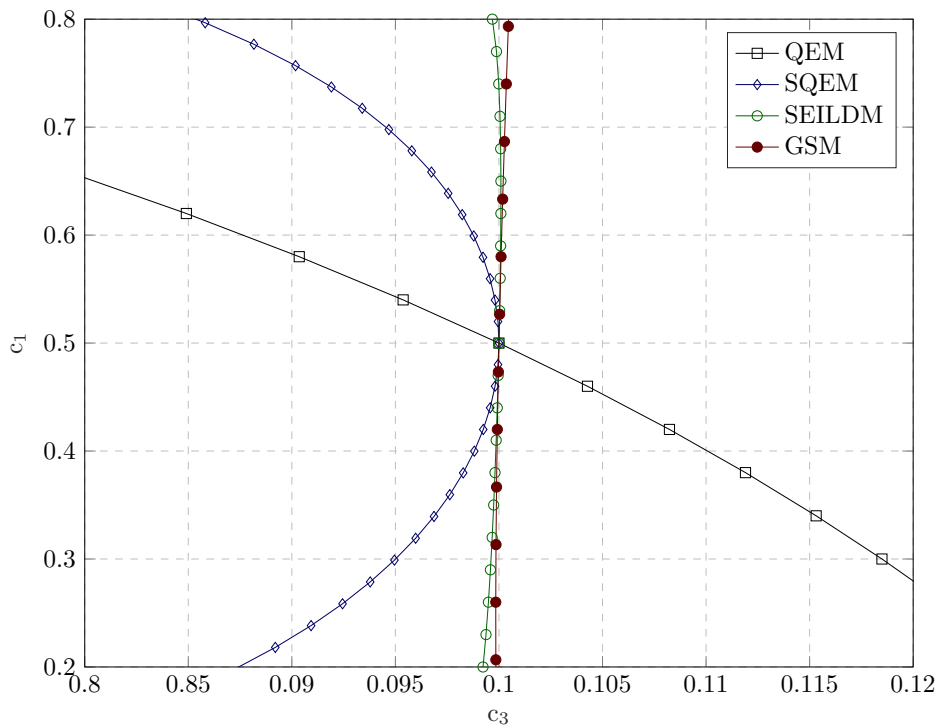
\begin{figure}[htb] \label{fig:compare_approx} 
	\centering
	\resizebox{0.97\textwidth}{!}{
%
%
\begin{tikzpicture}

\begin{axis}[%
width=4.521in,
height=3.538in,
at={(0.758in,0.509in)},
scale only axis,
xmin=0.08,
xmax=0.12,
xlabel style={font=\color{white!15!black}},
xlabel={$\text{c}_\text{3}$},
xticklabels={0,0,0.8,0.85,0.09,0.095,0.1,0.105,0.11,0.115,0.12},
ymin=0.2,
ymax=0.8,
ylabel style={font=\color{white!15!black}},
ylabel={$\text{c}_\text{1}$},
axis background/.style={fill=white},
xmajorgrids,
ymajorgrids,
grid style={dashed},
legend style={legend cell align=left, align=left, draw=white!15!black}
]
\addplot [color=black, mark=square, mark options={solid, black}]
  table[row sep=crcr]{%
0.131534156157351	0.1\\
0.129250468164976	0.14\\
0.126819067527067	0.18\\
0.124226193607375	0.22\\
0.121456473923901	0.26\\
0.118492709363268	0.3\\
0.115315630701689	0.34\\
0.111903623783509	0.38\\
0.108232421182164	0.42\\
0.104274759112199	0.46\\
0.1	0.5\\
0.0953737237579854	0.54\\
0.0903572948328543	0.58\\
0.0849074179894181	0.62\\
0.078975704079657	0.66\\
0.0725082782364625	0.7\\
0.0654454773288364	0.74\\
0.0577217009238293	0.78\\
0.0492654992519836	0.82\\
0.04	0.86\\
0.0298437881283575	0.9\\
};
\addlegendentry{QEM}

\addplot [color=black!60!blue, mark=diamond, mark options={solid, black!60!blue}]
  table[row sep=crcr]{%
0.0719591298869628	0.0989116222850207\\
0.076009924936455	0.118352152813081\\
0.079456317418561	0.138038184767466\\
0.0824208424604951	0.157901643077503\\
0.084991987248791	0.177893368195116\\
0.0872355274335033	0.197976818092498\\
0.0892015121055829	0.218124187719121\\
0.0909287544425899	0.238313914198561\\
0.0924478218876674	0.258529015617963\\
0.0937830865231299	0.278755951931594\\
0.0949541661390711	0.298983824367183\\
0.0959769582115049	0.319203800993609\\
0.0968643945095282	0.339408697494707\\
0.0976269992205914	0.359592667099672\\
0.0982733056359539	0.379750969091137\\
0.0988101686339276	0.399879795203374\\
0.099242998498104	0.419976139723276\\
0.0995759337021649	0.440037703498797\\
0.0998119647843978	0.460062825119779\\
0.0999530174492049	0.480050434750442\\
0.0999999999999999	0.500000027777778\\
0.0999528177351695	0.519911656813795\\
0.0998103547307816	0.539785941816233\\
0.0995704212375451	0.559624099312475\\
0.0992296624802741	0.579427993066514\\
0.0987834216628681	0.599200210187296\\
0.098225546040807	0.618944168866218\\
0.0975481194205164	0.638664266988602\\
0.0967410964315002	0.658366085315833\\
0.0957918018419575	0.678056665643357\\
0.0946842394272011	0.697744894762666\\
0.0933981248593795	0.717442041744789\\
0.091907507440152	0.73716252364436\\
0.0901787605580028	0.756925021912221\\
0.0881675696580327	0.776754155745537\\
0.0858142654818557	0.796683074732302\\
0.0830362998128469	0.816757638992863\\
0.0797155145803904	0.837043491899783\\
0.0756752872634754	0.857638757075847\\
0.0706363521144386	0.878698582158645\\
0.0641229723989941	0.900487265333892\\
};
\addlegendentry{SQEM}

\addplot [color=black!60!green, mark=o, mark options={solid, black!60!green}]
  table[row sep=crcr]{%
0.0984294167797731	0.11\\
0.0987964890343818	0.14\\
0.0990462909522835	0.17\\
0.0992305234675807	0.2\\
0.0993742222424479	0.23\\
0.0994909567909214	0.26\\
0.099588690784918	0.29\\
0.0996723740273406	0.32\\
0.0997452093834825	0.35\\
0.0998093199207488	0.38\\
0.0998661195013938	0.41\\
0.0999165249540503	0.44\\
0.0999610765938546	0.47\\
0.1	0.5\\
0.1	0.5\\
0.100033223903955	0.53\\
0.100060357186027	0.56\\
0.100080617065877	0.59\\
0.100092685552085	0.62\\
0.100094444063752	0.65\\
0.100082478764336	0.68\\
0.100051113726784	0.71\\
0.0999903732250467	0.74\\
0.099881211108427	0.77\\
0.0996825865705432	0.8\\
0.099288171451515	0.83\\
0.0983237972948938	0.86\\
0.0945480926730073	0.89\\
};
\addlegendentry{SEILDM}

\addplot [color=black!60!red, mark=*, mark options={solid, black!60!red}]
  table[row sep=crcr]{%
0.09990380859375	0.1\\
0.0998681488037109	0.153333333333333\\
0.0998491821289062	0.206666666666667\\
0.0998460388183594	0.26\\
0.0998577117919922	0.313333333333333\\
0.0998831787109375	0.366666666666667\\
0.0999213409423828	0.42\\
0.0999711151123047	0.473333333333333\\
0.100031372070312	0.526666666666667\\
0.100101028442383	0.58\\
0.100179046630859	0.633333333333333\\
0.100264373779297	0.686666666666667\\
0.100356079101563	0.74\\
0.100453231811523	0.793333333333333\\
0.100554977416992	0.846666666666667\\
0.100660537719727	0.9\\
};
\addlegendentry{GSM}

\end{axis}
\end{tikzpicture}%
}
	\caption{Illustration of different SIM approximation methods applied to the test system defined by \eqref{eq:c_MM_chemical} and \eqref{eq:c_MM} respectively, similar to Figure 4(a) in \cite{chiavazzo2007comparison}. A SIM is given by the vertical line $c_3 = 0.1$.}
\end{figure} the  following SIM approximation methods are applied in \cite{chiavazzo2007comparison}:
\begin{itemize}
	\item Quasi-Equilibrium-Manifold (QEM)
	\item Spectral-Quasi-Equilibrium-Manifold (SQEM)
	\item Symmetric-Entropic-Intrinsic-Low-Dimensional-Manifold (SEILDM)
\end{itemize}
Considering the conservation law and choosing a specific set of constants for that mechanism (see \cite{chiavazzo2007comparison} for more details), we receive the two-dimensional system
\begin{equation}\label{eq:c_MM}
\begin{cases}
\dot{c}_3 &= c_3^2-2.1c_3+0.2\\
\dot{c}_1 &= 0.5c_3+c_1c_3-0.2c_1
\end{cases}
\end{equation}
We apply our Geodesic Stretching Method (GSM) to this system. In this case, we have a one-dimensional SIM - according to the GSPT-definition. The SIM is the vertical line $c_3 = 0.1$ with equilibrium point $(c_1 = 0.5, c_3 = 0.1)$. In Figure \ref{fig:compare_approx} you can see the corresponding approximations plotted in the $(c_3-c_1)$-plane. Applying the GSM, $c_1$ is our reaction-progress variable.

We can see that the approximation error of the QEM and the SQEM is comparatively large. Both SEILDM and GSM provide significantly better approximations with similar deviations from the SIM. One of two candidates provided by the solution of the ILDM equation coincides exactly with the SIM which is unsurprising, since the SIM is linear and the approximation error of the ILDM is proportional to the curvature of the SIM \cite{kaper2002asymptotic}.

\subsubsection*{General applicability of Geodesic Stretching}
The previous two test models are comparatively simple in the following manner: Both are
\begin{itemize}
	\item[(a)] two-dimensional with a one-dimensional SIM,
	\item[(b)] slow-fast systems,
	\item[(c)] formulated on the (linear) phase-space $\R^2$.
\end{itemize} 
We briefly point out how this approach can be applied to systems that are not limited to the properties (a) - (c) in our test examples, it can be used for systems of every (finite) dimension $n$ with a SIM of an arbitrary dimension $1\le  m < n$ (addressing point (a)). Definitions \ref{def:pure_space_division} and \ref{def: tangVSorth_stret} already cover how this more general case is tackled: Consider the subspace $T_{p,s}^\text{tan}(U)$ of the whole tangent space that is tangent to a potential SIM $U$ and calculate the maximum stretching rate for each tangent vector of this subspace (denoted as $\varTheta_p^{\text{tan}}(U)$). This generalization is still covariant since every tangent vector is a tensor with a coordinate-independent meaning.

In case of a general system $\dot{x} = f(x)$, the object we want to approximate is a NAIM, as the notion of a SIM is tied to slow-fast systems. The geodesic stretching approach is formulated for a system of the form $\dot{x}= f(x)$ and approximates normal attractiveness. Hence, there is no limitation to slow-fast systems and Restriction (b) is addressed by construction.

Concerning point (c), imagine a model containing adiabatic constraints $g(x) =0$ where $g:\R^n \to \R^\ell$ with $1 \le \ell <n$. Instead of operating on an open subset of $\R^n$, the dynamic is now formulated on a manifold. Conveniently, all exploited notions and tensors used to define this geodesic stretching approach are formulated on Riemann manifolds. We only need to choose a parent coordinate system and define the metric on this manifold according to \cref{def:metric_tensor}. This approach can be applied directly, once this is done.

\section{ Differential Geometric Interpretation of the Flow Curvature Method}\label{sec:geo_ginoux}
Sections 2 and 3 cover a new covariant approach on manifolds to approximate NAIM's core property by utilizing a differential-geometry setting. In this section, we aim to investigate a potential covariant reformulation of an existing geometry-based approach to SIMs: The flow-curvature method (FCM) by Ginoux \cite{ginoux2006differential}. We begin by briefly stating its main properties:

\subsection{Flow Curvature Method in a nutshell} \label{sec:FCM_nutshell}
The foundation of this ansatz are higher curvatures of trajectories in the phase space $\R^n$. In case of an $n$-dimensional dynamical system  $ \dot{x} = f (x) $  with $(n-1)$-dimensional SIM, its FCM-approximation is defined by the union of all points $p$ with vanishing $n$-th curvature of the trajectory. It can be shown that the FCM is capable of approximating a SIM to order $n$ \cite{ginoux2008slow}. The former criterion is satisfied if and only if
\[
\Psi(x):=\det \left( \frac{\d}{\d t}x(t), \ldots , \frac{\d^{(n)}}{\d t^{(n)}}x(t)  \right)\bigg\vert_{x=p} = 0.
\]
Its solution is called flow curvature manifold  
which is not flow-invariant, as long as $\frac{\d}{\d t} J_f(x(t))$ does not vanish. This can directly be seen by calculating $\frac{\d}{\d t} \left(\Psi(x(t)) \right).$ 
	
On the other hand, this manifold is also non-invariant regarding coordinate transformations. Away from each fixed point, we can locally transform the system into a constant system, e.g. $ \dot{y} = g(y)\equiv c, $ with $  y = \Phi^{(c)}(x), $ for a given system $\dot{x} = f(x)$.
 In these new coordinates, every point satisfies the flow curvature criterion. This argument also proves that the FCM is not covariant. The value of $\Psi$ - and crucially - whether or not $\Psi$ vanishes depends on the coordinate chart, as $\Phi^{(c)}$ shows.

An evident way of finding a covariant reformulation of the FCM is to take steps similar to the approach used to develop the geodesic stretching method. Step one: Translate the existing method into a manifold-based setting. Step two: Modify the method in a sensible way, such that the scalar value (here: $\Psi$) is the evaluation of a tensor and can be expressed in any appropriate coordinate chart. We implement the first step in the following subsections.
\subsection{Flow derivatives as covariant derivatives in euclidean space}
Consider the manifold $(\M, g^e) = (\R^n,\langle \cdot , \cdot \rangle)$ where $\langle \cdot , \cdot \rangle = g^{e}$ represents the euclidean inner product at each point $p \in \M$. 
Let $\partial_{j,p} \in T_p \M$ indicates the tangent vector in the direction of 
the $j$-th coordinate. Let $\d x_{j,p}$ represent the dual basis on each point $p$, we receive
\[   
g_p^{(e)} = \sum\limits_{j = 1}^n \d x_{j,p} \otimes \d x_{j,p}\qquad \forall p \in \R^n.
\]
The christoffel symbols of the Levi-Civita connection all vanish. 
\begin{lemma} \label{lem: cvd_der_gin}
	Let $(\M, g^{e})$ be given as above and $f:\R^n \to \R^n$ be sufficiently smooth. Let $\nabla= \nabla^{e}$ be the Levi-Civita connection preserving $g^e$. Suppose $(\d x / \d t)= f(x)$. Let $h:\R^n \to \R^n $ be continuously differentiable. Then the flow derivative of $h$ coincides with the covariant derivative in the direction $f(x)$: 
	\[
	\nabla_{f(p)} h = \frac{\d}{\d t} h(x(t))\big\vert_p \qquad \forall p \in \R^n.
	\]
\end{lemma}
\begin{proof}
	Direct calculation, see Appendix. 
\end{proof}
We define a matrix column-wise consisting of the first $n$ flow derivatives  
\[
M(p):= \left[  \frac{\d}{\d t}x(t), \ldots , \frac{\d^{(n)}}{\d t^{(n)}}x(t) \right]\bigg\vert_{x=p}.
\]
Let the successive covariant derivative be denoted by
\[
\nabla^{(\ell)}_\beta \alpha := \underbrace{\nabla_{\beta} \ldots \nabla_{\beta}}_{\ell -\text{times}}(\alpha)\qquad \text{and} \qquad \nabla_{\beta}^{(0)}(\alpha):= \alpha
\]
for sufficiently smooth vector fields $\beta$ and $\alpha$ on $\R^n$ and $\ell \in \mathbb{N}$. \cref{lem: cvd_der_gin} implies
\begin{corollary}
	Let $(\R^n,g^e)$ and $\nabla$ be defined as above. Suppose that $\dot{x}=f(x)$, then we get
	\[
	\frac{\d^{(k+1)} }{\d t^{(k+1)}} (x(t)) = \nabla^{(k)}_f (f) \qquad \forall k \in \mathbb{N},
	\]
	implying that we can rewrite $M(p)$ as
	\begin{equation}\label{eq:cov_flow_mat}
	M(p) =  \left[ \nabla_{f}^{(0)}f, \ldots , \nabla_{f}^{(n-1)}f \right]\bigg\vert_{x=p}.
	\end{equation}
\end{corollary}
\begin{proof}
	Using the previous lemma iteratively, we get
	\begin{align*}
	\frac{\d^{(k+1)}}{\d t^{(k+1)}}(x(t)) &= \frac{\d}{\d t}\left( \frac{\d^{(k)}}{\d t^{(k)}} x(t) \right) = \nabla_{f(p)}\left( \frac{\d^{(k)}}{\d t^{(k)}} x(t) \right) = ... \\
	&= \underbrace{\nabla_{f(p)}... \nabla_{f}}_{k-\text{times}} \left( \frac{\d}{\d t}(x(t))  \right) = \nabla_{f}^{(k)}(f)
	\end{align*}
\end{proof}
In the FCM criterion use \cref{eq:cov_flow_mat} to reformulate the function $\Phi$:
\[
0 = \det(M(p))= \det \left( \nabla_{f}^{(0)} f, \ldots, \nabla_{f}^{(n-1)} f\right)\bigg\vert_{x=p}.
\]
\subsection{Flow Curvature Function as Gramian Determinant}

\begin{definition}
	Let $v_1,\ldots v_n$ be vector fields on $\R^n$. The Gramian matrix $G_p:(T_p \R^n)^n \to \R^{n\times n}$ and Gramian determinant $D_p: (T_p \R^n)^n \to \R$ are defined by
	\begin{align*}
	(v_{1,p}, \ldots, v_{n,p}) &\mapsto G_p(v_{1,p}, \ldots, v_{n,p}) := \left(g_p^e(v_{i,p},v_{j,p}) \right)_{i,j}\\
	(v_{1,p}, \ldots, v_{n,p}) &\mapsto D_p(v_{1,p}, \ldots, v_{n,p}) := \sqrt{\det\left(\left(g_p^e(v_{i,p},v_{j,p}) \right)_{i,j}\right)}
	\end{align*}
	respectively.
\end{definition}
By definition, both the gramian matrix and determinant are coordinate independent for every metric $g$. In case of the euclidean metric $g^e$, a direct calculation shows
\begin{equation}
G_p \left( \frac{\d }{\d t}x(t)\bigg \vert_{x=p}, \ldots ,  \frac{\d^{(n)} }{\d t^{(n)}}x(t)\bigg \vert_{x=p} \right) = M(p)^T M(p).
\end{equation}
Using the multiplicativity of the determinant, we conclude that the definition criterion for the FCM can be written in the following coordinate independent way:
\[
\Phi(p) = 0 \Leftrightarrow D_p(\nabla_{f}^{(0)}f, \ldots, \nabla_{f}^{(n-1)}f) = 0
\]

The last transformation finishes the embedding of the FCM into the Riemann geometry framework. All translated quantities yield are well-defined within this field. Unfortunately, the utilized notions are not tensors - and as mentioned in the \cref{sec:FCM_nutshell} - a modification is now needed. One might be guided by the geometric interpretation of the FCM: Inspect the highest curvature of a solution trajectory. Find a way to translate this interpretation into a differential geometry framework to define the appropriate tensors.
\section{Summary, Conclusion and Outlook}
The aim of this work was to present a route to formulate tensorial representation of normal attractiveness in multiple time-scale systems, enabling tensorial approximations for slow invariant manifolds for slow-fast systems. For this purpose we introduce a general differential geometric setting in \Cref{sec:Diffgeo_chapter}. 
 We exploit the notion of \textit{intrinsic curvature}  to reformulate the stretching-based analysis. This covariant formulation on manifolds makes this work a novelty in this context.
We exemplarily apply the resulting approach to the Davis-Skodje system and the Michaels-Menten model in \Cref{sec:test_stret}.
 In \Cref{sec:geo_ginoux} we also reformulate the flow curvature method by expressing its utilized flow derivatives by covariant ones. Our ideas might be useful as a general guideline towards finding tensorial reformulations of established SIM methods.
The authors share the opinion that the field of differential geometry is an appropriate frame as to express essential SIM quantities intrinsically as e.g. shown in \cite{heiter2018towards}.
\section*{Acknowledgments}
The authors thank the Klaus-Tschira foundation (project 00.003.2019) for financial funding, as well as Marcus Heitel and Jörn Dietrich for discussions on this topic.
\section*{Appendix}
\subsection*{Appendix A: Christoffel symbols and Proof of Theorem}
\subsubsection*{Inverse metric tensor:} Let $p \in \M$ be given and $g$ be a metric tensor - a symmetric positive bilinear form - on $T_p \M$ be given. The mapping \[
\phi: T_p \M \to T_p^{\prime} \M, \qquad v_p \mapsto g_p(v_p, \cdot) \in  T_p^{\prime} \M \qquad \forall v_p \in T_p \M
\]
is an isomorphism. There exists an unique symmetric, positive, bilinear mapping
\[
\tilde{g}: T_p^{\prime} \M \times T_p^{\prime} \M \to \R
\]
such that the mapping 
\[
\psi:  T_p^{\star} \M \to T_p \M, \qquad d_p \to \iota(\tilde{g}(d_p, \cdot)), \qquad \forall d_p \in T_p^{\prime}M
\]
is the inverse of $\phi$. Here, $\iota$ represents the natural identification of $T_p^{\prime}M$ and its bidual space $T_p^{\prime \prime}M$.
The component matrix of the so-called inverse metric tensor is denoted by $(g^{ij}_p)_{i,j}$ and satisfies the equality
\[
(g^{ij}_p)_{i,j} (g_{ij,p})_{i,j} = \text{Id}_{n+1}.
\]
Inverting the component matrix $(g_{ij,p})_{i,j}$ of $g_p$ w.r.t. the basis $\{ \partial_{1,p}, \ldots, \partial_{n+1,p} \}$ from chapter 2 leads to
\begin{equation}\label{eq:metric_inv}
g_p^{ij} = 
\begin{pmatrix}
\text{Id}_n + f(x_p)f(x_p)^T & f(x_p)\\
f(x_p)^T & 1		
\end{pmatrix}			
\end{equation}
\subsubsection*{Calculation of Christoffel symbols:}
For the sake of simplicity, the dependence on $p\in \M$ is left out in the following calculations.
Because the chosen connection is the Levi-Civita connection, the Christoffel symbols can be calculated directly by the formula
\begin{equation}\label{eq:christ_formula}
\Gamma_{ij}^k = \frac{1}{2} g^{k \ell} \left( \partial_i g_{j \ell} + \partial_j g_{i \ell} - \partial_{\ell} g_{ij} \right)\qquad \forall k,i,j \in \{ 1,\ldots,n+1 \}. 
\end{equation}
Let $k \in \{1,\ldots,n\}$ and let $f_{\ell} = f_{\ell}(x_p)$ indicate the $\ell$-th entry of the vector $f(x_p)$. Plugging the components $g_{ij}$ from equation \cref{eq:metr_comps} into \cref{eq:christ_formula} yields
\[
(\Gamma^k_{ij})_{i,j} = 
\left[
\begin{array}{c|c}
\begin{array}{c}
-\frac{f_k}{2}\left( \frac{\delta f_i}{\delta x_j} + \frac{\delta f_j}{\delta x_i} \right) \\
(i,j) \in \{1,\ldots, n\}^2  
\end{array}&
\begin{array}{c}
\left(\frac{f_k}{2} \sum_{\mu} f_{\mu} \left( \frac{\delta f_j}{\delta x_{\mu}} + \frac{\delta f_{\mu}}{\delta x_j}\right) \right) + \frac{1}{2} \left( \frac{\delta f_j}{\delta x_{k}} - \frac{\delta f_{k}}{\delta x_j}\right) \\
j = 1,\ldots, n
\end{array} \\
\hline
* & (-f_k)f^T J_f f - \sum_{\mu} f_{\mu}\frac{\delta f_{\mu}}{\delta x_k}
\end{array}
\right]
\] 
where the entries marked by a $*$ are determined by the symmetry $\Gamma^k_{ij} = \Gamma^k_{ji}$. In case $k=n+1$, the components we receive are:
\[
(\Gamma^{n+1}_{ij})_{i,j} = 
\left[
\begin{array}{c|c}
\begin{array}{c}
-\frac{1}{2}\left( \frac{\delta f_i}{\delta x_j} + \frac{\delta f_j}{\delta x_i} \right) \\
(i,j) \in \{1,\ldots, n\}^2  
\end{array} & 
\begin{array}{c}
\frac{1}{2} \sum_{\mu} f_{\mu} \left( \frac{\delta f_j}{\delta x_{\mu}} + \frac{\delta f_{\mu}}{\delta x_j}\right)\\j = 1,\ldots, n
\end{array} \\
\hline
* & -f^T J_f f 
\end{array}
\right].
\]
\subsubsection*{Proof of Theorem \cref{thm:geodetization}}
\begin{proof}
	Let $\gamma: (-\varepsilon,\varepsilon) \to \M$, $t\mapsto \gamma(t)$ with $\gamma(0) = [x_0,\tau_0]$ be a solution trajectory of the extended system. The first and second derivative of $\gamma$ w.r.t. t are given by
	\[
	\frac{\d \gamma }{\d t}(0) = \begin{pmatrix}
	f(x_0)\\
	1
	\end{pmatrix}\qquad  \frac{\d^2 \gamma }{\d t^2}(0) = \begin{pmatrix}
	J_f(x_0) f(x_0) \\0
	\end{pmatrix}.
	\]
	Let $k \in \{1, \ldots, n\}$. We calculate \begin{align*}
	&\left(\frac{\d}{\d t}( x(t), \tau(t)) \right) (\Gamma_{ij}^k(\gamma(t)))_{i,j} \frac{\d}{\d t} 
	\begin{pmatrix}
	x(t)\\
	\tau(t)
	\end{pmatrix}\\
	= &
	(f^T,1) (\Gamma^k_{ij})_{i,j} \begin{pmatrix}
	f \\1
	\end{pmatrix}\\
	=& f^T \left(\begin{array}{c}
	-\frac{f_k}{2}\left( \frac{\delta f_i}{\delta x_j} + \frac{\delta f_j}{\delta x_i} \right) \\
	(i,j) \in \{1,\ldots, n\}^2   
	\end{array}\right)f+ (-f_k)f^T J_f f - \sum_{\mu} f_{\mu}\frac{\delta f_{\mu}}{\delta x_k} \\
	+& 2 f^T  \left(\left(\frac{f_k}{2} \sum_{\mu} f_{\mu} \left( \frac{\delta f_j}{\delta x_{\mu}} + \frac{\delta f_{\mu}}{\delta x_j}\right) \right) + \frac{1}{2} \left( \frac{\delta f_j}{\delta x_{k}}- \frac{\delta f_{k}}{\delta x_j}\right)\right)_{j=1,\ldots, n+1} \\
	=& -2 f_k( f^T J_f f) + 2 f_k(f^T J_f f)- \sum_{\mu} f_{\mu}\frac{\delta f_{\mu}}{\delta x_k}+ \sum_{\mu} f_{\mu}\left(\frac{\delta f_{\mu}}{\delta x_k}
	- \frac{\delta f_{k}}{\delta x_{\mu}}\right)\\
	=& - \sum_{\mu} f_{\mu}\frac{\delta f_{k}}{\delta x_{\mu}} = -\frac{\d^2 \gamma^{(k)} }{\d t^2}
	\end{align*}
	In the case that $k=n+1$, we receive
	\begin{align*}
	&\left(\frac{\d}{\d t}( x(t), \tau(t)) \right) (\Gamma_{ij}^k(\gamma(t)))_{i,j} \frac{\d}{\d t} 
	\begin{pmatrix}
	x(t)\\
	\tau(t)
	\end{pmatrix}\\
	&= f^T \left( \begin{array}{c}
	-\frac{1}{2}\left( \frac{\delta f_i}{\delta x_j} + \frac{\delta f_j}{\delta x_i} \right) \\
	(i,j) \in \{1,\ldots, n\}^2  
	\end{array}  \right) f - f^T J_f f  \\
	&+ 2 f^T \left( \frac{1}{2} \sum_{\mu} f_{\mu} \left( \frac{\delta f_j}{\delta x_{\mu}} + \frac{\delta f_{\mu}}{\delta x_j}\right) \right)_{j = 1, \ldots, n} \\
	&= -2 f^T J_f f + 2 f^T J_f f = 0.
	\end{align*}
	Insertion of the identities from above into the geodesic equation proves the Theorem.
\end{proof}

\subsection*{Appendix B: Proof of \cref{lem: cvd_der_gin}}
\begin{proof}
	
	The euclidean metric $g^{e}$ satisfies $g_p^e(\partial_{i,p},\partial_{j,p}) = \delta_{ij}$ for all tuples $(i,j)\in \{1, \ldots, n\}^2$, implying 
	\[
	\frac{\partial g_{ij}^e}{\partial x_k} = 0\quad \forall (i,j,k) \in \{1, \ldots,n\}^3 \Rightarrow \Gamma_{ij}^k = 0 \quad \forall (i,j,k) \in \{1, \ldots,n\}^3.
	\]
	Thus, the covariant derivatives of the base vector fields $\nabla_{\partial_{i,p}} \partial_{j}$ vanish. Let $h \in T \R^n$ be a smooth vector field with components $h_k(x)$. We receive
	\begin{align*}
	\nabla_{\partial_{i,p}}h &= \nabla_{\partial_{i,p}}\left(\sum_{k=1}^{n} h_k(x) \partial_{k,p} \right) = \sum_{k=1}^{n}\nabla_{\partial_{i,p}}\left( h_k(x)  \partial_{k,p} \right) \\
	&= \sum_{k=1}^{n}\left( h_k(x) \nabla_{\partial_{i,p}} \partial_{k,p} + \frac{h_k}{\partial x_i} \partial_{k,p} \right) =  \sum_{k=1}^{n} \frac{h_k}{\partial x_i} \partial_{k,p} \qquad \forall i = 1, \ldots, n.
	\end{align*}
	By linearity we conclude that 
	\[
	\nabla_{f(p)}h = \sum_{i=1}^n f_i(p) \sum_{k=1}^n \frac{\partial h_k}{\partial x_{i}} \partial_{k,p}.
	\]
	The component vector of the flow derivatives is calculated by
	\[
	\frac{\d h(x(t))}{\d t} \bigg \vert_p = J_h(p) f(p) = \sum_{i=1}^n f_i(p) \sum_{k=1}^n \frac{\partial h_k}{\partial x_{i}}
	\]
	which proves the Lemma.
\end{proof}

\bibliographystyle{plain}
\bibliography{references_JDSGT}

\begin{thebibliography}{10}

\bibitem{adrover2007structure}
Alessandra Adrover, Francesco Creta, Stefano Cerbelli, Massimiliano Valorani,
  and Massimiliano Giona.
\newblock The structure of slow invariant manifolds and their bifurcational
  routes in chemical kinetic models.
\newblock {\em Computers \& Chemical Engineering}, 31(11):1456--1474, 2007.

\bibitem{adrover2007stretching}
Alessandra Adrover, Francesco Creta, Massimiliano Giona, and Massimiliano
  Valorani.
\newblock Stretching-based diagnostics and reduction of chemical kinetic models
  with diffusion.
\newblock {\em Journal of Computational Physics}, 225(2):1442--1471, 2007.

\bibitem{arnol1963proof}
Vladimir~I Arnol'd.
\newblock Proof of a theorem of an kolmogorov on the invariance of
  quasi-periodic motions under small perturbations of the hamiltonian.
\newblock {\em Russian Mathematical Surveys}, 18(5):9, 1963.

\bibitem{bodenstein1913theorie}
Max Bodenstein.
\newblock Eine theorie der photochemischen reaktionsgeschwindigkeiten.
\newblock {\em Zeitschrift f{\"u}r physikalische Chemie}, 85(1):329--397, 1913.

\bibitem{Chapman1913}
David~Leonard Chapman and Leo~Kingsley Underhill.
\newblock The interaction of chlorine and hydrogen. the influence of mass.
\newblock {\em Journal of the Chemical Society, Transactions}, 103:496--508,
  1913.

\bibitem{chiavazzo2007comparison}
Eliodoro Chiavazzo, Alexander~N Gorban, and Iliya~V Karlin.
\newblock Comparison of invariant manifolds for model reduction in chemical
  kinetics.
\newblock {\em Commun. Comput. Phys}, 2(5):964--992, 2007.

\bibitem{davis1999geometric}
Michael~J Davis and Rex~T Skodje.
\newblock Geometric investigation of low-dimensional manifolds in systems
  approaching equilibrium.
\newblock {\em The Journal of chemical physics}, 111(3):859--874, 1999.

\bibitem{eldering2013normally}
Jaap Eldering.
\newblock {\em Normally hyperbolic invariant manifolds: the noncompact case},
  volume~2.
\newblock Springer, 2013.

\bibitem{fenichel1974asymptotic}
Neil Fenichel.
\newblock Asymptotic stability with rate conditions.
\newblock {\em Indiana University Mathematics Journal}, 23(12):1109--1137,
  1974.

\bibitem{fenichel1977asymptotic}
Neil Fenichel.
\newblock Asymptotic stability with rate conditions, ii.
\newblock {\em Indiana University Mathematics Journal}, 26(1):81--93, 1977.

\bibitem{fenichel1979geometric}
Neil Fenichel.
\newblock Geometric singular perturbation theory for ordinary differential
  equations.
\newblock {\em Journal of differential equations}, 31(1):53--98, 1979.

\bibitem{fenichel1971persistence}
Neil Fenichel and JK~Moser.
\newblock Persistence and smoothness of invariant manifolds for flows.
\newblock {\em Indiana University Mathematics Journal}, 21(3):193--226, 1971.

\bibitem{gear2005projecting}
C~William Gear, Tasso~J Kaper, Ioannis~G Kevrekidis, and Antonios Zagaris.
\newblock Projecting to a slow manifold: Singularly perturbed systems and
  legacy codes.
\newblock {\em SIAM Journal on Applied Dynamical Systems}, 4(3):711--732, 2005.

\bibitem{ginoux2006differential}
Jean-Marc Ginoux and Bruno Rossetto.
\newblock Differential geometry and mechanics: applications to chaotic
  dynamical systems.
\newblock {\em International Journal of Bifurcation and Chaos},
  16(04):887--910, 2006.

\bibitem{ginoux2008slow}
Jean-Marc Ginoux, Bruno Rossetto, and Leon~O Chua.
\newblock Slow invariant manifolds as curvature of the flow of dynamical
  systems.
\newblock {\em International Journal of Bifurcation and Chaos},
  18(11):3409--3430, 2008.

\bibitem{hadamard1901iteration}
Jacques Hadamard.
\newblock Sur l’it{\'e}ration et les solutions asymptotiques des
  {\'e}quations diff{\'e}rentielles.
\newblock {\em Bull. Soc. Math. France}, 29:224--228, 1901.

\bibitem{heiter2018towards}
Pascal Heiter and Dirk Lebiedz.
\newblock Towards differential geometric characterization of slow invariant
  manifolds in extended phase space: Sectional curvature and flow invariance.
\newblock {\em SIAM Journal on Applied Dynamical Systems}, 17(1):732--753,
  2018.

\bibitem{hirsch2006invariant}
Morris~W Hirsch, Charles~Chapman Pugh, and Michael Shub.
\newblock {\em Invariant manifolds}, volume 583.
\newblock Springer, 2006.

\bibitem{kaper2002asymptotic}
Hans~G Kaper and Tasso~J Kaper.
\newblock Asymptotic analysis of two reduction methods for systems of chemical
  reactions.
\newblock {\em Physica D: Nonlinear Phenomena}, 165(1-2):66--93, 2002.

\bibitem{kaper2015geometry}
Hans~G Kaper, Tasso~J Kaper, and Antonios Zagaris.
\newblock Geometry of the computational singular perturbation method.
\newblock {\em Mathematical modelling of natural phenomena}, 10(3):16--30,
  2015.

\bibitem{kuehn2015multiple}
Christian Kuehn.
\newblock {\em Multiple time scale dynamics}, volume 191.
\newblock Springer, 2015.

\bibitem{lam1989understanding}
SH~Lam and DA~Goussis.
\newblock Understanding complex chemical kinetics with computational singular
  perturbation.
\newblock In {\em Symposium (International) on Combustion}, volume~22, pages
  931--941. Elsevier, 1989.

\bibitem{lam1994csp}
SH~Lam and DA~Goussis.
\newblock The csp method for simplifying kinetics.
\newblock {\em International journal of chemical kinetics}, 26(4):461--486,
  1994.

\bibitem{lebiedz2010entropy}
Dirk Lebiedz.
\newblock Entropy-related extremum principles for model reduction of
  dissipative dynamical systems.
\newblock {\em Entropy}, 12(4):706--719, 2010.

\bibitem{lebiedz2010minimal}
Dirk Lebiedz, Volkmar Reinhardt, and Jochen Siehr.
\newblock Minimal curvature trajectories: Riemannian geometry concepts for slow
  manifold computation in chemical kinetics.
\newblock {\em Journal of Computational Physics}, 229(18):6512--6533, 2010.

\bibitem{lebiedz2011geometric}
Dirk Lebiedz, Volkmar Reinhardt, Jochen Siehr, and Jonas Unger.
\newblock Geometric criteria for model reduction in chemical kinetics via
  optimization of trajectories.
\newblock In {\em Coping with Complexity: Model Reduction and Data Analysis},
  pages 241--252. Springer, 2011.

\bibitem{lebiedz2013continuation}
Dirk Lebiedz and Jochen Siehr.
\newblock A continuation method for the efficient solution of parametric
  optimization problems in kinetic model reduction.
\newblock {\em SIAM Journal on Scientific Computing}, 35(3):A1584--A1603, 2013.

\bibitem{lebiedz2014optimization}
Dirk Lebiedz and Jochen Siehr.
\newblock An optimization approach to kinetic model reduction for combustion
  chemistry.
\newblock {\em Flow, Turbulence and Combustion}, 92(4):885--902, 2014.

\bibitem{lebiedz2011variational}
Dirk Lebiedz, Jochen Siehr, and Jonas Unger.
\newblock A variational principle for computing slow invariant manifolds in
  dissipative dynamical systems.
\newblock {\em SIAM Journal on Scientific Computing}, 33(2):703--720, 2011.

\bibitem{Lebiedz2016}
Dirk Lebiedz and Jonas Unger.
\newblock On unifying concepts for trajectory-based slow invariant attracting
  manifold computation in kinetic multiscale models.
\newblock {\em Mathematical and Computer Modelling of Dynamical Systems},
  22(2):87--112, 2016.

\bibitem{lee2006riemannian}
John~M Lee.
\newblock {\em Riemannian manifolds: an introduction to curvature}, volume 176.
\newblock Springer Science \& Business Media, 2006.

\bibitem{lyapunov1992general}
Aleksandr~Mikhailovich Lyapunov.
\newblock The general problem of the stability of motion.
\newblock {\em International journal of control}, 55(3):531--534, 1992.

\bibitem{maas1992simplifying}
Ulrich Maas and Stephen~B Pope.
\newblock Simplifying chemical kinetics: intrinsic low-dimensional manifolds in
  composition space.
\newblock {\em Combustion and flame}, 88(3-4):239--264, 1992.

\bibitem{poincare1899methodes}
Henri Poincar{\'e}.
\newblock {\em Les m{\'e}thodes nouvelles de la m{\'e}canique c{\'e}leste},
  volume~3.
\newblock Gauthier-Villars, 1899.

\bibitem{rosa1996inertial}
Ricardo Rosa and Roger Temam.
\newblock Inertial manifolds and normal hyperbolicity.
\newblock {\em Acta Applicandae Mathematica}, 45(1):1--50, 1996.

\bibitem{Roussel2012}
Marc~R. Roussel.
\newblock Further studies of the functional equation truncation approximation.
\newblock {\em Canadian Applied Mathematics Quarterly}, 20(2):209--227, 2012.

\bibitem{roussel1991geometry}
Marc~R Roussel and Simon~J Fraser.
\newblock On the geometry of transient relaxation.
\newblock {\em The Journal of chemical physics}, 94(11):7106--7113, 1991.

\bibitem{temam2012infinite}
Roger Temam.
\newblock {\em Infinite-dimensional dynamical systems in mechanics and
  physics}, volume~68.
\newblock Springer Science \& Business Media, 2012.

\bibitem{valorani2009g}
Mauro Valorani and Samuel Paolucci.
\newblock The g-scheme: A framework for multi-scale adaptive model reduction.
\newblock {\em Journal of Computational Physics}, 228(13):4665--4701, 2009.

\bibitem{wiggins1994normally}
Stephen Wiggins.
\newblock {\em Normally hyperbolic invariant manifolds in dynamical systems},
  volume 105.
\newblock Springer Science \& Business Media, 1994.

\bibitem{wiggins2013normally}
Stephen Wiggins.
\newblock {\em Normally hyperbolic invariant manifolds in dynamical systems},
  volume 105.
\newblock Springer Science \& Business Media, 2013.

\bibitem{zagaris2004analysis}
Antonios Zagaris, Hans~G Kaper, and Tasso~J Kaper.
\newblock Analysis of the computational singular perturbation reduction method
  for chemical kinetics.
\newblock {\em Journal of Nonlinear Science}, 14(1):59--91, 2004.

\bibitem{zagaris2005two}
Antonios Zagaris, Hans~G Kaper, and Tasso~J Kaper.
\newblock Two perspectives on reduction of ordinary differential equations.
\newblock {\em Mathematische Nachrichten}, 278(12-13):1629--1642, 2005.

\end{thebibliography}
\end{document}